\documentclass{amsart}
\usepackage{amsmath, amssymb, amsthm, amsfonts, amsgen, stmaryrd}

\numberwithin{equation}{section}

\makeatletter
\def\rightharpoonupfill@{\arrowfill@\relbar\relbar\rightharpoonup}
\newcommand{\xrightharpoonup}[2][]{\ext@arrow 0359\rightharpoonupfill@{#1}{#2}}
\makeatother

\newcommand{\res}{\mathop{\hbox{\vrule height 7pt width .5pt depth 0pt
\vrule height .5pt width 6pt depth 0pt}}\nolimits}
\newcommand{\ds}{\displaystyle}

\newcommand{\A}{{\mathcal{A}}}
\newcommand{\C}{{\mathcal{C}}}
\newcommand{\E}{{\mathcal{E}}}
\newcommand{\I}{{\mathcal{I}}}
\newcommand{\J}{{\mathcal{J}}}
\newcommand{\G}{{\mathcal{G}}}
\newcommand{\M}{{\mathcal{M}}}
\newcommand{\F}{{\mathcal{F}}}
\newcommand{\LL}{{\mathcal{L}}}
\newcommand{\HH}{{\mathcal{H}}}

\newcommand{\Nb}{{\mathbb{N}}}

\newcommand{\Rb}{{\mathbb{R}}}

\newcommand{\medi}{- \hskip-0.9em \int}
\newcommand{\med}{- \hskip-1.0em \int}

\let\e=\varepsilon
\let\a=\alpha
\let\b=\beta
\let\d=\delta

\let\O=\Omega
\let\o=\omega
\let\G=\Gamma

\let\wto=\rightharpoonup

\newtheorem{thm}{Theorem}[section]

\newtheorem{rmk}[thm]{Remark}
\newtheorem{lemma}[thm]{Lemma}
\newtheorem{proposition}[thm]{Proposition}

\begin{document}

\thispagestyle{empty}
\title[Lower semicontinuity of quasiconvex bulk energies in SBV]
{Lower semicontinuity of quasiconvex bulk energies in SBV\\and
integral representation in dimension reduction}
\author{Jean-Fran\c cois Babadjian}

\address[J.-F. Babadjian]{SISSA, Via Beirut 2-4, 34014
Trieste, Italy} \email{babadjia@sissa.it}

\maketitle

\begin{abstract}
A result of Larsen concerning the structure of the approximate
gradient of certain sequences of functions with Bounded Variation is
used to present a short proof of Ambrosio's lower semicontinuity
theorem for quasiconvex bulk energies in $SBV$. It enables to
generalize to the $SBV$ setting the decomposition lemma for scaled
gradients in dimension reduction and also to show that, from the
point of view of bulk energies, $SBV$ dimensional reduction problems
can be reduced to analogue ones in the Sobolev spaces framework.\\

\keywords{\noindent {\sc Keywords:}} Dimension reduction,
$\Gamma$-convergence, functions of bounded variation, free
discontinuity problems, quasiconvexity, equi-integrability.\\

\subjclass{\noindent {\sc 2000 Mathematics Subject Classification:}
74K35, 49J45, 49Q20.}
\end{abstract}


\section{Introduction}

\noindent Since the pioneering work \cite{LDR}, the modelling of
thin films through dimensional reduction techniques and
$\G$-convergence analysis has become one of the main issues in the
field of Calculus of Variations. In the membrane theory framework in
nonlinear elasticity, the problem rests on the study of the (scaled)
elastic energy
$$\frac{1}{\e} \int_{\O_\e} W(\e)(y,\nabla v)\, dy$$
of such bodies. Here $\O_\e:=\o \times (-\e/2,\e/2)$, where $\o$ is
a bounded open subset of $\Rb^2$ and $\e>0$, stands for the
reference configuration of a nonlinear elastic thin film, $v:\O_\e
\to \Rb^3$ is the deformation field which maps the reference
configuration into a deformed configuration and $W(\e) : \O_\e
\times \Rb^{3 \times 3} \to [0,+\infty)$ is the stored energy
density of the body which is a Carath\'eodory function satisfying
uniform $p$-growth and $p$-coercivity conditions (with
$1<p<\infty$). From a mathematical point of view, the previous
energy is well defined provided $v$ is a Sobolev function in
$W^{1,p}(\O_\e;\Rb^3)$.

To study the limit problem as the thickness $\e \to 0$, it will be
useful to recast the energy functional over the varying set $\O_\e$
into a functional with a fixed domain of integration $\O:=\o \times
(-1/2,1/2)$. To this end, denoting by $x_\a:=(x_1,x_2)$ the in-plane
variable, we set $u(x_\a,x_3):=v(x_\a,\e\, x_3)$ so that, after the
(now standard) change of variables
$$x_\a=y_\a, \quad x_3=\frac{y_3}{\e},$$
we are equivalently led to study the following rescaled functional
\begin{equation}\label{Sobolev}
\int_{\O} W_\e\left(x,\nabla_\a u \Big|\frac{1}{\e}\nabla_3 u\right)
dx, \quad u \in W^{1,p}(\O;\Rb^3),
\end{equation}
where $W_\e:\O \times \Rb^{3 \times 3} \to [0,+\infty)$ is the
rescaled stored energy density expressed in the new variables and
defined by $W_\e(x_\a,x_3,\xi):=W(\e)(x_\a,\e\, x_3,\xi)$. From now
on, $\nabla_\a$ (resp. $\nabla_3$) will stand for the (approximate)
gradient with respect to $x_\a$ (resp. $x_3$), $\xi=(\xi_\a|\xi_3)$
for some matrix $\xi \in \Rb^{3 \times 3}$ and $z=(z_\a|z_3)$ for
some vector $z \in \Rb^3$. Thus in view of the $p$-growth of the
energy, it is important to understand the structure of what we call
the scaled gradient of $u$, {\it i.e.}
\begin{equation}\label{scaledgrad}
\left(\nabla_\a u \Big|\frac{1}{\e}\nabla_3 u \right).
\end{equation}
In particular, if $\{u_\e\} \subset W^{1,p}(\O;\Rb^3)$ is a
minimizing sequence uniformly bounded in energy, up to a
subsequence, there always exist $u \in W^{1,p}(\O;\Rb^3)$ such that
$D_3u=0$ in the sense of distributions and $b \in L^p(\O;\Rb^3)$
such that $u_\e \wto u$ in $W^{1,p}(\O;\Rb^3)$ and $(1/\e)\nabla_3
u_\e \wto b$ in $L^p(\O;\Rb^3)$. The limit function $u$ is nothing
but the deformation of the mid-plane while $b$ is called the {\it
Cosserat vector}. It seems thus natural to expect a limit model
depending on the pair $(u,b)$. Unfortunately, this is still out of
reach and we refer to \cite{FFL} for a more detailed discussion on
the subject. However, in \cite{BFMbis} (see also \cite{babF}) a
simplified model has been considered taking into account the {\it
bending moment} $\overline b \in L^p(\o;\Rb^3)$, {\it i.e.} the
average in the transverse direction
$x_3$ of $b$, instead of the full Cosserat vector field.\\

In the framework of fracture mechanics, one usually adds a surface
energy term, penalizing the presence of the crack. The simplest case
consists in just penalizing its area leading to the so-called
Griffith's surface energy.  Thus, for a given crack, one should
study the energy given by the competing sum of the bulk and the
surface energies. Such fracture mechanics problems belong (among
others) to the class of free discontinuity problems, that is
variational problems where the unknown is not only a function, but a
pair set/function. Based on the idea that the deformation may be
discontinuous across the crack, it is convenient to study the weak
formulation, replacing the crack by the jump set of the deformation
and leading to a variational problem stated in the space of
(Special) Functions with Bounded Variation. Now the energy in which
we are interested is
$$\frac{1}{\e} \int_{\O_\e} W(\e)(y,\nabla v)\, dy + \frac{1}{\e} \HH^2(S_v), \quad v \in SBV^p(\O_\e;\Rb^3),$$
where $\nabla v$ is intended as the approximate gradient of $v$,
$S_v$ is the jump set of $v$ and $\HH^2$ stands for the
$2$-dimensional Hausdorff measure. Writing as before this energy in
the rescaled variables yields to
\begin{equation}\label{SBV}
\int_{\O} W_\e\left(x,\nabla_\a u \Big|\frac{1}{\e}\nabla_3
u\right)dx + \int_{S_u} \left|\left(\left( \nu_u\right)_\a
\Big|\frac{1}{\e}\left(\nu_u\right)_3 \right) \right|d\HH^2, \quad u
\in SBV^p(\O;\Rb^3)
\end{equation}
where $\nu_u$ is the generalized normal to $S_u$ and
(\ref{scaledgrad}) is now referred as the approximate scaled
gradient of $u$.\\

The aim of this paper is to study the connections between
variational problems (\ref{Sobolev}) and (\ref{SBV}), possibly
taking into account the presence of the bending moment vector field.
To this end, we will use as main ingredient Theorem \ref{bab} which
extends the Decomposition Lemma for scaled gradients (see
\cite[Theorem~1.1]{BF} or \cite[Theorem~3.1]{BZ}) to the $SBV$
setting. It states that any $SBV$-sequence with bounded rescaled
bulk energy and whose derivative's singular part behaves
asymptotically well, can be energetically  replaced, up to a set of
vanishing Lebesgue measure, by a sequence of Lipschitz maps whose
scaled gradient is $p$-equi-integrable. Thus it reduces the free
discontinuity problem to a usual dimensional reduction one in the
framework of Sobolev spaces. This result is nothing but a rescaled
version of \cite[Lemma~2.1]{L} (see also Theorem \ref{Larsen}
below). Using this structure theorem, we are able to show two
integral representation theorems in $SBV$ (Theorems
\ref{gammaconvbend} and \ref{gammaconv}) which say that, up to a
subsequence, the functional (\ref{SBV}) $\G$-converges (in an
appropriate topology) to a functional of the same kind, {\it i.e.}
the sum of a bulk and a surface energy. Moreover, the surface energy
is still of Griffith's type while the bulk energy is exactly the
same than that obtained in the analogue Sobolev spaces analysis. The
main importance of these representation theorems relies on the fact
that results on dimension reduction in Sobolev spaces can now be
extended to $SBV$ (see \cite{babF,BB1,BB2,BFMbis,LDR}).

Note that an integral representation result for dimensional
reduction problems in $SBV$ already exists (see
\cite[Theorem~2.1]{BraFo}). Even if this reference may seem more
general from the point of view of the hypothesis, it does not
contain as special case our results because the authors made
strongly use of the fact that their surface energy had to grow
linearly with respect to the deformation jump. This assumption was
essential in order to get compactness in $BV(\O;\Rb^3)$ of
minimizing sequences. However, they suggested a way to remove that
constraint by singular perturbation \cite[Remark 2.2]{BraFo}. In our
study we use a direct argument based on a trick introduced in
\cite{FF} and which was already used in \cite{babadjian} in the
framework of dimensional reduction. It consists in defining an
artificial functional exactly as we usually do for the
$\G$-$\liminf$, except that we impose the minimizing sequences to be
uniformly bounded in $L^\infty(\O;\Rb^3)$. Thanks to a truncation
argument (see Lemma \ref{bd}) we show that it actually coincides
with the $\G$-$\liminf$ for deformations $u \in L^\infty(\O;\Rb^3)$
and the advantage is that now, minimizing sequences turn out to be
relatively compact in $SBV(\O;\Rb^3)$ thanks to Ambrosio's
Compactness Theorem. We refer to \cite{babadjian} for a deeper
insight on that subject.

To close this introduction, we wish to stress that in this paper, we
are mostly interested in representation of effective bulk energies
arising in 3D-2D dimensional reduction problems stated in $SBV$. For
this reason we will consider a large class for such bulk energies
while surface energies will be restricted to the simplified case of
a Griffith's type one. However we are convinced that the results
presented here could be generalized to a larger class of surface
energies.\\

The overall plan of the paper is as follows: after recalling some
useful notations in section \ref{notation} and in order to show the
technique in a more transparent way, we present in section \ref{amb}
a short proof of Ambrosio's lower semicontinuity result for
quasiconvex integrands using \cite[Lemma~2.1]{L}. Then in section
\ref{lars} we prove our main tool, Theorem \ref{bab}, thanks to a
slicing argument together with \cite[Lemma~2.1]{L}. To reach our
goal, we need to prove a general integral representation for the
$\G$-limit of (\ref{Sobolev}) in $W^{1,p}(\O;\Rb^3) \times
L^p(\o;\Rb^3)$ as a function of the deformation and the bending
moment. This is the purpose of Theorem \ref{bending} in section
\ref{IRS} which contains as particular cases
\cite[Theorem~3.1]{BFMbis} (with $W_\e(x,\xi)=W(\xi)$) and
\cite[Theorem~3.4]{babF} (with $W_\e(x,\xi)=W(x,\xi)$). In section
\ref{gammaconvergence}, we refine the analysis of section \ref{amb}
adding the difficulties of dimension reduction. From the integral
representation in Sobolev spaces, Theorem \ref{bending}, we deduce
an analogue result in $SBV$, Theorem \ref{gammaconvbend}, which says
that the $\G$-limit of (\ref{SBV}) in $BV(\O;\Rb^3) \times
L^p(\o;\Rb^3)$ has also an integral representation and that the bulk
energy is exactly the same one than that obtained in the $W^{1,p}$
analysis. This will be achieved thanks to Theorem \ref{bab} and a
blow-up method which enables to reduce the problem to affine
deformations and constant bending moments. Finally we deduce a
similar result in section \ref{7} without the presence of the
bending moment.

\section{Notations and preliminaries}\label{notation}

\noindent If $\O \subset \Rb^N$ is an open set, we consider the
Lebesgue spaces $L^p(\O;\Rb^d)$ and the Sobolev spaces
$W^{1,p}(\O;\Rb^d)$ in the usual way. When needed, we will precise
what topology the space $L^p(\O;\Rb^d)$ will be endowed. In
particular we will denote by $L^p_s(\O;\Rb^d)$ (resp.
$L^p_w(\O;\Rb^d)$) the space $L^p(\O;\Rb^d)$ endowed with the strong
(resp. weak) topology. Strong convergence will always be denoted by
$\to$ while weak (resp. weak*) convergence will be denoted by $\wto$
(resp. $\xrightharpoonup[]{*}$).

We denote by $\M(\O;\Rb^d)$ the space of vector valued finite Radon
measures. If $\mu \in \M(\O;\Rb^d)$ and $E$ is a Borel subset of
$\O$, we will write $\mu \res \,  E$ for the restriction of $\mu$ to
$E$ that is, for every Borel subset $F$ of $\O$, $\mu\, \res \,  E\,
(F)=\mu(E \cap F)$. The Lebesgue measure in $\Rb^N$ will be denoted
by $\LL^N$ while $\HH^{N-1}$ is the $(N-1)$-dimensional Hausdorff
measure. We will denote by $B$ the unit ball of $\Rb^N$ and by
$\o_N:=\LL^N(B)$ its Lebesgue measure. If $x_0 \in \Rb^N$ and
$\rho>0$, $B(x_0,\rho):=x_0 +\rho\, B$ is the ball centered at $x_0$
with radius $\rho$. The notation $\medi_A$ stands for the average
$\LL^N(A)^{-1}\int_A$.

The space of Functions of Bounded Variation is denoted by
$BV(\O;\Rb^d)$ and we refer to \cite{AFP} for standard theory of
$BV$ functions. We recall here few facts: if $u \in BV(\O;\Rb^d)$
then its distributional derivative $Du \in \M(\O;\Rb^{d \times N})$
and thanks to Lebesgue's Decomposition Theorem, we can write $Du=D^a
u + D^s u$, where $D^au$ and $D^su$ stand for, respectively, the
absolutely continuous and singular part of $Du$ with respect to the
Lebesgue measure $\LL^N$. Let $S_u$ be the complementary of Lebesgue
points of $u$. We say that $u$ is a Special Function of Bounded
Variation, and we write $u \in SBV(\O;\Rb^d)$, if
$$Du=\nabla u\, \LL^N + (u^+ - u^-) \otimes \nu_u \,
\HH^{N-1}\res \,  S_u$$ where $\nabla u$ is the approximate gradient
of $u$, $\nu_u$ is the generalized normal to $S_u$ and $u^\pm$ are
the traces of $u$ on both sides of $S_u$. If $E \subset \O$, we say
that $E$ has finite perimeter in $\O$ provided $\chi_E \in SBV(\O)$.
We denote by $\partial^* E$ (resp. $\partial_* E$) the reduced
(resp. essential) boundary of $E$. When $p>1$, we define
$$SBV^p(\O;\Rb^d) := \Big\{u \in SBV(\O;\Rb^d): \nabla u \in L^p(\O;\Rb^{d \times N}) \text{ and }
\HH^{N-1}(S_u \cap \O) <+ \infty \Big\}.$$ We say that a sequence
$\{u_n\} \subset SBV^p(\O;\Rb^d)$ converges weakly to some $u \in
SBV^p(\O;\Rb^d)$, and we write $u_n \wto u$ in $SBV^p(\O;\Rb^d)$, if
$$\left\{
\begin{array}{l}
u_n \to u \text{ in } L^1(\O;\Rb^d),\\[0.2cm]
\nabla u_n  \rightharpoonup \nabla u \text{ in }L^p(\O;\Rb^{d \times N}),\\[0.1cm]
(u_n^+-u_n^-) \otimes \nu_{u_n} \HH^{N-1}\res \,  S_{u_n}
\xrightharpoonup[]{*} (u^+-u^-) \otimes \nu_{u} \HH^{N-1}\res \, S_u
\text{ in }\M(\O;\Rb^{d \times N}).
\end{array}
\right.$$

If $\O:=\o \times I$, where $\o$ is a bounded open subset of $\Rb^2$
and $I:=(-1/2,1/2)$, we will identify the spaces $L^p(\o;\Rb^3)$,
$W^{1,p}(\o;\Rb^3)$ or $SBV^p(\o;\Rb^3)$ with the space of functions
$v \in L^p(\O;\Rb^3)$, $W^{1,p}(\O;\Rb^3)$ or $SBV^p(\O;\Rb^3)$)
such that $D_3 v=0$ in the sense of distributions.

By $\A(\o)$ we mean the family all open subsets of $\o$ while
$\mathcal R(\o)$ stands for the countable subfamily of $\A(\o)$
obtained by taking all finite unions of open cubes contained in
$\o$, centered at rational points and with rational edge length.

In the sequel, we will denote by $Q':=(-1/2,1/2)^2$ the unit cube of
$\Rb^2$ and by $Q'(x_0,\rho):=x_0+\rho \,Q'$ the cube centered at
$x_0 \in \Rb^2$ and side length $\rho>0$. Similarly $B':=\{x_\a \in
\Rb^2: |x_\a| < 1\}$ stands for the unit ball in $\Rb^2$ and
$B'(x_0,\rho):=x_0+\rho\, B'$ denotes the ball of $\Rb^2$ centered
at $x_0 \in \Rb^2$ and of radius $\rho>0$.

\section{Lower semicontinuity of quasiconvex bulk energies in
$SBV$}\label{amb}

\noindent This section is devoted to give a short proof of
Ambrosio's lower semicontinuity result for quasiconvex bulk energies
in $SBV$ using the following theorem proved in \cite[Lemma~2.1]{L}.

\begin{thm}\label{Larsen}
Let $\O \subset \Rb^N$ be a bounded open set with Lipschitz boundary
and let $\{u_n\} \subset BV(\O;\Rb^d)$ be such that
$$\left\{
\begin{array}{l}
\ds \sup_{n \in \Nb} \|u_n\|_{BV(\O;\Rb^d)} < +\infty,\\[0.4cm]
\ds \sup_{n \in \Nb} \|\nabla u_n\|_{L^p(\O;\Rb^{d \times N})} <
+\infty \quad \text{for some }p>1,\\[0.3cm]
\ds |D^s u_n|(\O) \to 0.
\end{array}
\right.$$ Then there exists a subsequence $\{n_k\} \nearrow +\infty$
and a sequence $\{w_k\} \subset W^{1,\infty}(\O;\Rb^d)$ such that
$$\left\{
\begin{array}{l}
\ds \sup_{k \in \Nb} \|w_k\|_{W^{1,p}(\O;\Rb^d)} < +\infty,\\[0.4cm]
\ds \{|\nabla w_k|^p\} \text{ is equi-integrable},\\[0.2cm]
\ds \LL^N(\{w_k \neq u_{n_k}\} \cup \{\nabla w_k \neq \nabla
u_{n_k}\}) \to 0.
\end{array}
\right.$$
\end{thm}

This theorem is nothing but the $BV$ counterpart of the
Decomposition Lemma, \cite[Lemma~1.2]{FMP}, in Sobolev spaces. We
now use the previous result to give a short proof of Ambrosio's
lower semicontinuity result for quasi-convex bulk energies in $SBV$
(see \cite[Theorem~4.3]{A2} or \cite[Proposition~5.29]{AFP}). This
will enable us to emphasize the techniques used in this paper,
occulting the difficulties of dimension reduction. The same kind of
arguments will be used in section \ref{gammaconvergence} to prove
the lower bound of Theorem \ref{gammaconvbend}.

\begin{thm}\label{Ambrosio}
Let $\O$ be bounded open subset of $\Rb^N$ and $f:\O \times \Rb^d
\times \Rb^{d \times N} \to [0,+\infty)$ be a Carath\'eodory
function satisfying
\begin{equation}\label{pgf}
c |\xi|^p \leq f(x,s,\xi) \leq a(x) + \psi(|s|)(1+|\xi|^p) \text{
for all }(s,\xi) \in \Rb^d \times \Rb^{d \times N} \text{ and a.e.
}x\in \O,
\end{equation}
for some $p>1$, $c>0$, $a \in L^1(\O)$ and some increasing function
$\psi:[0,+\infty) \to [0,+\infty)$. If $\xi \mapsto f(x,s,\xi)$ is
quasiconvex for every $s \in \Rb^d$ and a.e. $x\in \O$, then
$$\liminf_{n \to +\infty}\int_\O f(x,u_n,\nabla u_n)\, dx \geq \int_\O f(x,u,\nabla u)\,
dx$$ for any sequence $\{u_n\} \subset SBV(\O;\Rb^d)$ converging in
$L^1(\O;\Rb^d)$ to $u \in SBV(\O;\Rb^d)$ and satisfying $\sup_n
\HH^{N-1}(S_{u_n}) <+\infty$.
\end{thm}

\begin{proof}
The proof is divided into three steps. We first apply the blow-up
method to reduce the study to an affine limit function. Then we
prove that the resulting sequence can be modified, without
increasing too much the energy, into another one uniformly bounded
in $L^\infty$. Finally we apply Theorem \ref{Larsen} to replace this
last
sequence of $SBV$ functions by a sequence of Sobolev functions.\\

{\bf Step 1.} Up to a subsequence, there is no loss of generality to
assume the existence of nonnegative and finite Radon measures
$\lambda$ and $\mu \in \M(\O)$ such that $f(\cdot,u_n,\nabla u_n)\,
\LL^N \xrightharpoonup[]{*} \lambda$ and $\HH^{N-1}\res \,  S_{u_n}
\xrightharpoonup[]{*} \mu$ in $\M(\O)$. To prove Theorem
\ref{Ambrosio} it is enough to check that
$$\lambda(\O) \geq \int_\O f(x,u,\nabla u)\,dx$$
and thanks to Lebesgue's Differentiation Theorem, it suffices to
show that
$$\frac{d\lambda}{d\LL^N}(x_0) \geq f(x_0,u(x_0),\nabla u(x_0))$$
for $\LL^N$-a.e. $x_0 \in \O$. Select $x_0 \in \O$ such that
\begin{itemize}
\item[(a)] $x_0$ is a Lebesgue point of $u$ and $a$ and  a point of approximate
differentiability of $u$;
\item[(b)] The Radon-Nikod\'ym derivative of $\lambda$ with respect to $\LL^N$
exists and is finite;
\item[(c)] the following limit exists and
\begin{equation}\label{itemc}\lim_{\rho \to 0^+} \frac{\mu(B(x_0,\rho))}{\o_{N-1}
\rho^{N-1}}=0;
\end{equation}
\item[(d)] for any sequence $\{\rho_i\} \searrow 0^+$ there exists a
subsequence $\{\rho_{i(k)}\}$ and a $\LL^N$-negligible set $E
\subset B$ such that \begin{equation}\label{fk}\lim_{k \to
+\infty}f(x_0+\rho_{i(k)}y,u(x_0)+\rho_{i(k)}s,\xi)=f(x_0,u(x_0),\xi)\end{equation}
locally uniformly in $\Rb^d \times \Rb^{d \times N}$ for any $y \in
B \setminus E$.
\end{itemize}
Note that $\LL^N$-a.e. points $x_0$ in $\O$ satisfy these
properties. Items (a) and (b) are immediate while item (d) is a
consequence of \cite[Lemma~5.38]{AFP}. Concerning item (c), we
remark that, setting
$$\Theta(x):=\limsup_{\rho \to 0^+}
\frac{\mu(B(x,\rho))}{\o_{N-1} \rho^{N-1}},$$ then $\{\Theta >0\} =
\bigcup_{h=1}^{+\infty}\{\Theta \geq 1/h\}$ and using
\cite[Theorem~2.56]{AFP}, we get that $\HH^{N-1}(\{\Theta \geq
1/h\}) \leq h \mu(\{\Theta \geq 1/h\}) < +\infty$. Thus
$\LL^N(\{\Theta \geq 1/h\})=0$ and consequently
$\LL^N(\{\Theta>0\})=0$.

Consider a sequence $\{\rho_k\} \searrow 0^+$ such that
$0<\rho_k<1$, $\mu(\partial B(x_0,\rho_k))=\lambda(\partial
B(x_0,\rho_k))=0$ for every $k\in \Nb$ and (\ref{fk}) holds with
$\rho_k$ in place of $\rho_{i(k)}$. Then
\begin{eqnarray}\label{1413}
\frac{d\lambda}{d\LL^N}(x_0) & = & \lim_{k \to
+\infty}\frac{\lambda(B(x_0,\rho_k))}{\o_N\rho_k^N}\nonumber\\
& = & \lim_{k \to +\infty} \lim_{n \to
+\infty}\frac{1}{\o_N\rho_k^N}\int_{B(x_0,\rho_k)}f(x,u_n,\nabla
u_n)\, dx\nonumber\\
& = & \lim_{k \to +\infty} \lim_{n \to
+\infty}\frac{1}{\o_N}\int_{B}f(x_0+\rho_k y,u(x_0)+\rho_k
u_{n,k},\nabla u_{n,k})\, dy
\end{eqnarray}
where we set $u_{n,k}(y)=[u_n(x_0+\rho_k\, y) - u(x_0)]/\rho_k$.
Since $x_0$ is a point of approximate differentiability of $u$, it
follows that
\begin{equation}\label{1414}
\lim_{k \to +\infty} \lim_{n \to +\infty}\|u_{n,k}
-w_0\|_{L^1(B;\Rb^d)} =0
\end{equation}
where $w_0(y):=\nabla u(x_0)\, y$. Moreover, by (\ref{itemc}) we get
that
\begin{eqnarray}\label{1415}
\lim_{k \to +\infty} \lim_{n \to +\infty}\HH^{N-1}(S_{u_{n,k}} \cap
B) & = & \lim_{k \to +\infty} \lim_{n \to
+\infty}\frac{\HH^{N-1}(S_{u_n}
\cap B(x_0,\rho_k))}{\rho_k^{N-1}}\nonumber\\
& =  & \lim_{k \to +\infty} \frac{\mu(B(x_0,\rho_k))}{\rho_k^{N-1}}=
0.
\end{eqnarray}
From (\ref{1413}), (\ref{1414}) and (\ref{1415}), one can find a
sequence $n(k) \nearrow +\infty$ such that, setting
$v_k:=u_{n(k),k}$, then $v_k \to w_0$ in $L^1(B;\Rb^d)$,
$\HH^{N-1}(S_{v_k} \cap B) \to 0$ and
\begin{equation}\label{v_k}
\frac{d\lambda}{d\LL^N}(x_0)=\lim_{k \to +\infty}
\frac{1}{\o_N}\int_{B}f(x_0+\rho_k y,u(x_0)+\rho_k v_k,\nabla v_k)\,
dy.
\end{equation}
From now on, all the integrals will be restricted to the unit ball
$B$.\\

{\bf Step 2.} We now use the same truncation argument than in the
proof of \cite[Proposition~5.37]{AFP}. Define $\hat
v_k:=(\sqrt{1+|v_k-w_0|^2}-2)^+$ so that by Theorem 3.96 and
Proposition 3.64 (c) in \cite{AFP}, $\hat v_k \in SBV(B)$, $|\nabla
\hat v_k| \leq |\nabla v_k - \nabla w_0|$ $\LL^N$-a.e. in $B$ and
$S_{\hat v_k} \subset S_{v_k}$. According the Coarea Formula in $BV$
\cite[Theorem~3.40]{AFP}, we have that \begin{eqnarray*}\int_0^1
\HH^{N-1}\big(\partial^*\{\hat v_k > t\} \cap (B \setminus S_{\hat
v_k})\big)\, dt & \leq & |D\hat v_k|(B  \setminus S_{\hat v_k}) =
\int_B |\nabla \hat v_k|\, dx\\
& \leq & \int_{B \cap \{|v_k-w_0|>\sqrt 3\}} |\nabla v_k-\nabla
w_0|\, dx\end{eqnarray*} where we have used the fact that $\nabla
\hat v_k=0$ $\LL^N$-a.e. in $B \cap \{|v_k-w_0|\leq \sqrt 3\}$. From
(\ref{v_k}) and the $p$-coercivity condition (\ref{pgf}), the
sequence $\{\nabla v_k\}$ is uniformly bounded in $L^p(B;\Rb^{d
\times N})$ and since $p>1$, it is equi-integrable. Using the fact
that $\LL^N(B \cap \{|v_k-w_0|>\sqrt 3\}) \to 0$ we obtain that the
right hand side of the previous relation tends to zero as $k \to
+\infty$. Consequently, one can find $t_k \in (0,1)$ such that
$A_k:=\{\hat v_k>t_k\}$ has finite perimeter in $B$ and
\begin{equation}\label{Sbk}
\lim_{k \to +\infty}\HH^{N-1}(B \cap \partial^*A_k \setminus S_{\hat
v_k})=0.
\end{equation}
Define $\tilde v_k:=v_k\chi_{B \setminus A_k}+w_0\chi_{B \cap A_k}$
so that $\tilde v_k \to w_0$ in $L^1(B;\Rb^d)$. As $|\hat v_k| \leq
t_k <1$ in $B \setminus A_k$ it follows that $|v_k-w_0| \leq
2\sqrt{2}$ in $B \setminus A_k$ and thus
\begin{equation}\label{bli}
\|\tilde v_k\|_{L^\infty(B;\Rb^d)} \leq \|v_k\|_{L^\infty(B
\setminus A_k;\Rb^d)} + \|w_0\|_{L^\infty(B;\Rb^d)}\leq M
\end{equation} where
$M>0$ is independent of $k$. Denoting by $v_k^-$ the exterior trace
of $v_k$ on $\partial^* A_k \cap B$ oriented by the inner normal of
$A_k$, Remark 3.85 in \cite{AFP} implies that $|v_k^-(x)|\leq M$ for
$\HH^{N-1}$-a.e. $x \in \partial^* A_k \cap B$ and thus
$$ \int_{\partial^* A_k \cap B}|v_k^-|\, d\HH^{N-1} \leq  M \HH^{N-1}(\partial^*A_k \cap B) <+\infty$$
so that \cite[Theorem~3.84]{AFP} ensures that $\tilde v_k \in
SBV(B;\Rb^d)$. Since $S_{\tilde v_k} \subset S_{v_k} \cup \partial_*
A_k$, by (\ref{Sbk}) we get that \begin{eqnarray*} \HH^{N-1}(B \cap
S_{\tilde v_k}) & \leq & \HH^{N-1}(B \cap S_{v_k}) +  \HH^{N-1}(B
\cap
\partial_* A_k \setminus S_{v_k})\\
& \leq & \HH^{N-1}(B \cap S_{v_k}) + \HH^{N-1}(B \cap
\partial^* A_k \setminus S_{\hat v_k}) \to 0\end{eqnarray*} where we used the fact
that $S_{\hat v_k} \subset S_{v_k}$ and $\HH^{N-1}(B \cap \partial_*
A_k \setminus \partial^* A_k)=0$. Using the locality of approximate
gradients and the $p$-growth condition (\ref{pgf}), we get that
\begin{eqnarray*}
\int_{B}f(x_0+\rho_k y,u(x_0)+\rho_k \tilde v_k,\nabla \tilde v_k)\,
dy & = & \int_{B \setminus A_k }f(x_0+\rho_k y,u(x_0)+\rho_k
v_k,\nabla v_k)\, dy\\
&& + \int_{B \cap A_k }f(x_0+\rho_k y,u(x_0)+\rho_k w_0 ,\nabla
u(x_0)\, dy\\
& \leq & \int_{B}f(x_0+\rho_k y,u(x_0)+\rho_k v_k,\nabla v_k)\, dy\\
&&+ \int_{B \cap A_k}\left[a(x_0+\rho_k y)+\psi(|u(x_0)+\rho_k
w_0|)(1+|\nabla u(x_0)|^p) \right]dy.
\end{eqnarray*}
By the choice of $x_0$, the sequence $\{a(x_0+\rho_k \cdot)\}$ is
strongly converging in $L^1(B)$ to $a(x_0)$ and thus it is
equi-integrable. Hence as $\LL^N(A_k) \leq \LL^N(\{|v_k-w_0| \geq
\sqrt 3\}) \to 0$ we deduce that the second term on the right hand
side of the previous relation tends to zero as $k \to +\infty$ and
thanks to (\ref{v_k}) it follows that
\begin{equation}\label{barv_k}
\frac{d\lambda}{d\LL^N}(x_0)\geq \limsup_{k \to +\infty}
\frac{1}{\o_N}\int_{B}f(x_0+\rho_k y,u(x_0)+\rho_k \tilde v_k,\nabla
\tilde v_k)\, dy.
\end{equation}

{\bf Step 3.} By (\ref{bli}) we have that $|D^s \tilde v_k|(B) \leq
2 M \HH^{N-1}(S_{\tilde v_k} \cap B) \to 0$ while the $p$-coercivity
condition (\ref{pgf}) and item (b) imply that
$$\sup_{k \in \Nb} \|\nabla \tilde v_k\|_{L^p(B;\Rb^{d \times N})}
<+\infty.$$ Consequently the sequence $\{\tilde v_k\}$ fulfills the
assumptions of Theorem \ref{Larsen} so that considering a suitable
(not relabeled) subsequence, there exist a Lebesgue measurable set
$E_k \subset B$ and a sequence $\{w_k\} \subset
W^{1,\infty}(B;\Rb^d)$ such that $\{|\nabla w_k|^p\}$ is
equi-integrable, $w_k=\tilde v_k$ on $B \setminus E_k$ and
$\LL^N(E_k) \to 0$. From the proof of \cite[Lemma~2.1]{L}, it can
also be checked that $\sup_k \|w_k\|_{L^\infty(B;\Rb^d)} \leq M$. As
$$\int_B|w_k-w_0|\, dy \leq \int_{B \setminus E_k} |\tilde
v_k-w_0|\, dy + 2M \LL^N(E_k) \to 0$$ it follows that $w_k \to w_0$
in $L^1(B;\Rb^d)$ and defining the set $B_k^t:=\{x\in B: |\nabla
w_k(x)|\leq t\}$, relation (\ref{barv_k}) leads to
$$\frac{d\lambda}{d\LL^N}(x_0) \geq \limsup_{t \to +\infty}\limsup_{k
\to +\infty} \frac{1}{\o_N}\int_{B_k^t \setminus E_k}f(x_0+\rho_k
y,u(x_0)+\rho_k w_k,\nabla w_k)\, dy.$$  Using now (\ref{fk}) with
$\rho_{i(k)}=\rho_k$, we obtain that
$$\lim_{k \to +\infty}\int_{B_k^t \setminus E_k}|f(x_0+\rho_k y,u(x_0)+\rho_k
w_k,\nabla w_k)-f(x_0,u(x_0),\nabla w_k)|\, dy=0$$ for each $t>0$,
implying that
\begin{equation}\label{Yap}\frac{d\lambda}{d\LL^N}(x_0) \geq \limsup_{t \to +\infty}\limsup_{k \to +\infty}
\frac{1}{\o_N}\int_{B_k^t \setminus E_k}f(x_0,u(x_0),\nabla w_k)\,
dy.\end{equation} Since $\LL^N(E_k) \to 0$, according to the
$p$-growth condition (\ref{pgf}) we get that for every $t>0$,
\begin{equation}\label{Yap2}\int_{E_k \cap B_k^t}f(x_0,u(x_0),\nabla w_k)\,
dy \leq \big(a(x_0) + \psi(|u(x_0)|) (1+t^p) \big)\LL^N(E_k)
\xrightarrow[k \to +\infty]{} 0.\end{equation} On the other hand,
Chebyshev's Inequality ensures the existence of a constant $c>0$
(independent of $k$ and $t$) such that $\LL^N(B \setminus B_k^t)
\leq c/t^p \to 0$ as $t \to +\infty$, so that the equi-integrability
of $\{|\nabla w_k|^p\}$ yields to
\begin{equation}\label{Yap3}\sup_{k \in \Nb} \int_{B \setminus B_k^t} f(x_0,u(x_0),\nabla w_k)\, dy \leq \sup_{k \in \Nb}
\int_{B \setminus B_k^t} \big(a(x_0) + \psi(|u(x_0)|) (1+|\nabla
w_k|^p) \big)\, dy \xrightarrow[t \to +\infty]{} 0. \end{equation}
Gathering (\ref{Yap}), (\ref{Yap2}) and (\ref{Yap3}), we deduce that
$$\frac{d\lambda}{d\LL^N}(x_0) \geq \limsup_{k \to +\infty}
\frac{1}{\o_N}\int_B f(x_0,u(x_0),\nabla w_k)\, dy$$ and since $w_k
\wto w_0$ in $W^{1,p}(B;\Rb^d)$, we can apply
\cite[Theorem~II-4]{AF} to conclude that
$$\frac{d\lambda}{d\LL^N}(x_0) \geq f(x_0,u(x_0),\nabla u(x_0)).$$
\end{proof}

\section{Structure of approximate scaled gradients}\label{lars}

\noindent In this section we prove the following Theorem \ref{bab}
which is a similar result than Theorem \ref{Larsen} in the context
of dimension reduction. Note that it generalizes
\cite[Theorem~1.1]{BF} and \cite[Theorem~3.1]{BZ} (with obvious
changes for $n$D-$(n-k)$D dimensional reduction). Its proof relies
on a slicing argument similar to that used in
\cite[Theorem~3.1]{BZ}. It will be instrumental in section
\ref{gammaconvergence} to prove Theorem \ref{gammaconvbend} because
it will enable to replace $SBV$ minimizing sequences by Lipschitz
ones without increasing the energy.

From now on, $\O:=\o \times I$ where $\o$ is a bounded open subset
of $\Rb^2$ and $I:=(-1/2,1/2)$.

\begin{thm}\label{bab}
Assume that $\o$ has a Lipschitz boundary and $p>1$. Let $\{\e_n\}
\searrow 0^+$ and $\{u_n\} \subset SBV^p(\O;\Rb^3)$ be such that
\begin{equation}\label{bounds}
\sup_{n \in \Nb} \left\{ \|u_n\|_{L^\infty(\O;\Rb^3)}+ \int_\O
\left|\left(\nabla_\a u_n\Big|\frac{1}{\e_n} \nabla_3 u_n
\right)\right|^p dx \right\} < +\infty,
\end{equation}
\begin{equation}\label{meas} \int_{S_{u_n}} \left| \left( \big(
\nu_{u_n} \big)_\a \Big| \frac{1}{\e_n} \big( \nu_{u_n} \big)_3
\right) \right| \, d\mathcal H^2  \to 0\end{equation} and that $u_n
\rightharpoonup u \text{ in }SBV^p(\O;\Rb^3)$, $(1/\e_n)\nabla_3 u_n
\wto b$ in $L^p(\O;\Rb^3)$ for some $u \in W^{1,p}(\o;\Rb^3)$ and $b
\in L^p(\O;\Rb^3)$. Then there exist a subsequence $\{\e_{n_k}\}
\subset \{\e_n\}$ and a sequence $\{z_k\} \subset
W^{1,\infty}(\O;\Rb^3)$ such that $z_k \rightharpoonup u$ in
$W^{1,p}(\O;\Rb^3)$, $(1/\e_{n_k})\nabla_3 z_k \wto b$ in
$L^p(\O;\Rb^3)$, the sequence $\big\{\big|\big( \nabla_\a
z_k|\frac{1}{\e_{n_k}}\nabla_3 z_k\big) \big|^p \big\}$ is
equi-integrable and
$$\LL^3(\{z_k \neq u_{n_k}\} \cup \{\nabla z_k \neq \nabla u_{n_k}\}) \to
0.$$
\end{thm}

\begin{proof}
The proof is based on a slicing argument. We first come back to the
non rescaled cylinder $\O_{\e_n}=\o \times (-\e_n/2,\e_n/2)$ of
thickness $\e_n$ setting $v_n(x_\a,x_3):=u_n(x_\a,x_3/\e_n)$. It
follows that for each $n \in \Nb$, $v_n \in SBV^p(\O_{\e_n};\Rb^3)$
and changing variable in (\ref{bounds}) we get that
\begin{equation}\label{bounds2}
\sup_{n \in \Nb} \left\{ \|v_n\|_{L^\infty(\O_{\e_n};\Rb^3)}+
\frac{1}{\e_n} \int_{\O_{\e_n}} |\nabla v_n|^p\,  dx \right\} <
+\infty
\end{equation}
and \begin{equation}\label{meas2}\HH^2(S_{v_n}) =\e_n\int_{S_{u_n}}
\left| \left( \big( \nu_{u_n} \big)_\a \Big| \frac{1}{\e_n} \big(
\nu_{u_n} \big)_3 \right) \right| \, d\mathcal H^2.
\end{equation}
We now periodize the functions $v_n$ in the transverse direction
defining
$$\hat v_n(x_\a,x_3):=\left\{
\begin{array}{rcl}
v_n(x_\a,-\e_n-x_3) & \text{if} & -\e_n < x_3 \leq -\frac{\e_n}{2},\\[0.1cm]
v_n(x_\a,x_3) & \text{if} & -\frac{\e_n}{2}< x_3 < \frac{\e_n}{2},\\[0.1cm]
v_n(x_\a,\e_n-x_3) & \text{if} & \frac{\e_n}{2} \leq x_3 < \e_n.
\end{array}
\right.$$ Then $\hat v_n \in SBV^p\big(\o \times
(-\e_n,\e_n);\Rb^3\big)$ for each $n \in \Nb$ and from
(\ref{bounds2}) and (\ref{meas2}) it follows that
\begin{equation}\label{bounds3}
\sup_{n \in \Nb} \left\{ \|\hat v_n\|_{L^\infty(\o \times
(-\e_n,\e_n);\Rb^3)}+ \frac{1}{\e_n} \int_{\o \times (-\e_n,\e_n)}
|\nabla \hat v_n|^p\, dx\right\} < +\infty.
\end{equation}
and \begin{equation}\label{meas3}\HH^2(S_{\hat v_n})=
2\e_n\int_{S_{u_n}} \left| \left( \big( \nu_{u_n} \big)_\a \Big|
\frac{1}{\e_n} \big( \nu_{u_n} \big)_3 \right) \right| \, d\mathcal
H^2.
\end{equation}
We are now in a position to extend $\hat v_n$ by periodicity in the
$x_3$ direction. Note that we do not create any additional jump set
because periodicity ensures continuity at the interface of each
slice. Let
$$N_n:= \left\{
\begin{array}{ll}
\ds \frac{1}{4\e_n}-\frac{1}{2} & \text{if }\; \ds \frac{1}{4\e_n}-\frac{1}{2} \in \Nb,\\[3mm]
\ds \left[ \frac{1}{4\e_n}+\frac{1}{2} \right] & \text{otherwise}
\end{array}
\right.$$ where $[t]$ denotes the integer part of $t$. For every
$i\in \{-N_n, \ldots, N_n\}$, we set
$I_{i,n}:=\big((2i-1)\e_n,(2i+1)\e_n\big)$ and $\O_{i,n}:= \o \times
I_{i,n}$. Note that $N_n$ is the largest integer such that $\O \cap
\O_{i,n} \ne \emptyset$ for every $i\in \{-N_n, \ldots, N_n\}$. We
define the function $\tilde v_n$ on $\O(n):= \o \times
(-(2N_n+1)\e_n, (2N_n+1)\e_n)$ by extending $\hat v_n$ by
periodicity in the $x_3$ direction on $\O(n)$:
$$\tilde v_n(x_\a,x_3)=\hat v_n(x_\a,x_3-2i\e_n) \text{ if } x_3 \in I_{i,n}.$$
Since $\O \subset \O(n)$, $\tilde v_n \in SBV^p(\O;\Rb^3)$ and
thanks to (\ref{bounds3}) and the definition of $N_n$, we have that
\begin{equation}\label{bounds4}
\sup_{n \in \Nb} \left\{ \|\tilde v_n\|_{L^\infty(\O;\Rb^3)}+
\int_\O |\nabla \tilde v_n|^p\,  dx\right\} < +\infty
\end{equation}
while (\ref{meas3}) together with (\ref{meas}) imply that
\begin{eqnarray}\label{meas4}
\HH^2(S_{\tilde v_n}) \leq c\int_{S_{u_n}} \left| \left( \big(
\nu_{u_n} \big)_\a \Big| \frac{1}{\e_n} \big( \nu_{u_n} \big)_3
\right) \right| \, d\mathcal H^2 \to 0.
\end{eqnarray}
As a consequence of (\ref{bounds4}) and (\ref{meas4}), the sequence
$\{\tilde v_n\}$ fulfills the assumptions of Theorem \ref{Larsen}.
Hence there exist a subsequence $\{\e_{n_k}\} \subset \{\e_n\}$ and
a sequence $\{w_k\} \subset W^{1,\infty}(\O;\Rb^3)$ such that
$$\left\{
\begin{array}{l}
\ds \sup_{k \in \Nb} \|w_k\|_{W^{1,p}(\O;\Rb^3)} <
+\infty,\\[0.4cm]
\ds \{|\nabla w_k|^p\} \text{ is equi-integrable},\\[0.3cm]
\ds \LL^3(\{\tilde v_{n_k} \neq w_k\} \cup \{\nabla \tilde v_{n_k}
\neq \nabla w_k\}) \to 0.
\end{array}
\right.$$ From De La Vall\'ee Poussin's criterion, one can find an
increasing and continuous function $\vartheta : [0,+\infty) \to
[0,+\infty]$ such that $\vartheta(t)/t \to +\infty$ as $t \to
+\infty$ and
$$\sup_{k \in \Nb} \int_\O \vartheta (|\nabla w_k|^p)\,dx <
+\infty.$$ We claim that for at least half of the indexes $i \in
\{-N_{n_k}+1, \ldots, N_{n_k}-1\}$, there holds
\begin{equation}\label{index}
\frac{2N_{n_k}-1}{2} \int_{\O_{i,n_k}} \big[ \vartheta (|\nabla
w_k|^p) + |w_k|^p + |\nabla w_k|^p \big]\, dx \leq \int_{\O} \big[
\vartheta (|\nabla w_k|^p) + |w_k|^p + |\nabla w_k|^p \big]\, dx.
\end{equation}
If not, define $J_k$ to be the set of indexes $i \in  \{ -N_{n_k}+1,
\ldots, N_{n_k}-1\}$ such that (\ref{index}) does not hold. Then it
would imply that $\#(J_k)> (2N_{n_k}-1)/2$ and
\begin{eqnarray*}
\int_{\O} \big[ \vartheta (|\nabla w_k|^p) + |w_k|^p + |\nabla
w_k|^p \big]\, dx & \geq & \sum_{i \in J_k}\int_{\O_{i,n_k}} \big[
\vartheta (|\nabla w_k|^p) +
|w_k|^p + |\nabla w_k|^p \big]\, dx\\
& > & \frac{2}{2N_{n_k}-1} \#(J_k) \, \int_{\O} \big[ \vartheta
(|\nabla w_k|^p) + |w_k|^p + |\nabla w_k|^p \big]\, dx
\end{eqnarray*}
which is absurd. Similarly, one can show that for at least half of
the indexes satisfying (\ref{index}), we have that
\begin{equation}\label{index2}
\frac{2N_{n_k}-1}{4}\, \LL^3\res \,\O_{i,n_k} (\{\tilde v_{n_k} \neq
w_k\} \cup \{\nabla \tilde v_{n_k} \neq \nabla w_k\}) \leq
\LL^3(\{\tilde v_{n_k} \neq w_k\} \cup \{\nabla \tilde v_{n_k} \neq
\nabla w_k\}).
\end{equation}
Let $i_k \in \{ -N_{n_k}+1, \ldots, N_{n_k}-1\}$ be such that
(\ref{index}) and (\ref{index2}) hold at the same time. Define now
$z_k(x_\a,x_3):=w_k(x_\a,\e_{n_k} x_3+2\e_{n_k}\, i_k)$. Changing
variable in (\ref{index}) and (\ref{index2}) and using the
construction of $\tilde v_{n_k}$ from $u_{n_k}$ we get that
\begin{eqnarray*}&&\e_{n_k} \frac{2N_{n_k}-1}{2}
\int_{\O} \left[ \vartheta \left(\left|\left( \nabla_\a
z_k\Big|\frac{1}{\e_{n_k}}\nabla_3 z_k\right) \right|^p \right) +
|z_k|^p +\left|\left( \nabla_\a z_k\Big|\frac{1}{\e_{n_k}}\nabla_3
z_k\right) \right|^p
\right]\, dx\\
&&\hspace{2cm} \leq \int_{\O} \big[ \vartheta (|\nabla w_k|^p) +
|w_k|^p + |\nabla w_k|^p \big]\, dx
\end{eqnarray*} and
$$\e_{n_k} \frac{2N_{n_k}-1}{4}\, \LL^3\res \, \O\, (\{u_{n_k} \neq z_k\} \cup
\{\nabla u_{n_k} \neq \nabla z_k\}) \leq \LL^3(\{\tilde v_{n_k} \neq
w_k\} \cup \{\nabla \tilde v_{n_k} \neq \nabla w_k\}).$$ Since
$\e_{n_k}(2N_{n_k}-1) \geq 1/4$ for $k$ large enough, it follows
that
$$\left\{\begin{array}{l}
\ds \sup_{k \in \Nb} \int_{\O} \left[ \vartheta \left(\left|\left(
\nabla_\a z_k\Big|\frac{1}{\e_{n_k}}\nabla_3 z_k\right) \right|^p
\right) + |z_k|^p +\left|\left( \nabla_\a
z_k\Big|\frac{1}{\e_{n_k}}\nabla_3 z_k\right)
\right|^p \right]\, dx < +\infty,\\
\\
\ds \LL^3 (\{u_{n_k} \neq z_k\} \cup \{\nabla u_{n_k} \neq \nabla
z_k\}) \to 0 \end{array}\right.$$ and the equi-integrability of
$\big\{\big|\big( \nabla_\a z_k|\frac{1}{\e_{n_k}}\nabla_3 z_k\big)
\big|^p \big\}$ follows from De La Vall\'ee Poussin's criterion.

It remains to prove the weak convergence of $z_k$ and
$(1/\e_{n_k})\nabla_3 z_k$. Let $v \in L^{p'}(\O;\Rb^3)$ with
$1/p+1/p'=1$, then $$\int_\O (z_k - u)\cdot v\, dx
=\int_{\{z_k=u_{n_k}\}} (u_{n_k}-u)\cdot v \, dx + \int_{\{z_k \neq
u_{n_k}\}} (z_k-u)\cdot v \, dx.$$ As $\LL^3(\{z_k \neq
u_{n_k}\})\to 0$, it follows that $v \chi_{\{z_k = u_{n_k}\}} \to v$
in $L^{p'}(\O;\Rb^3)$. Then, using H\"older's Inequality, the fact
that $\{z_k\}$ is uniformly bounded in $L^p(\O;\Rb^3)$ and that
$u_{n_k} \wto u$ in $L^p(\O;\Rb^3)$, we obtain that
\begin{eqnarray*}
\lim_{k \to +\infty} \left|\int_\O (z_k - u)\cdot v\, dx \right| &
\leq  & \lim_{k \to +\infty} \left|\int_\O (u_{n_k}-u)\cdot
v\chi_{\{z_k = u_{n_k}\}} \, dx\right|\\
&&+ \lim_{k \to +\infty}\|z_k-u\|_{L^p(\O;\Rb^3)} \| v \chi_{\{z_k
\neq u_{n_k}\}}\|_{L^{p'}(\O;\Rb^3)} = 0.
\end{eqnarray*}
Similarly we may show that $\nabla z_k \wto \nabla u$ in
$L^p(\O;\Rb^{3 \times 3})$ and that $(1/\e_{n_k})\nabla_3 z_k \wto
b$ in $L^p(\O;\Rb^3)$.
\end{proof}

\section{Integral representation for dimension reduction
problems in Sobolev spaces involving the bending moment}\label{IRS}

\noindent Consider a Carath\'eodory function $W_\e:\O \times \Rb^{3
\times 3} \to [0,+\infty)$ satisfying uniform $p$-growth and
$p$-coercivity conditions: there exist $0<\b'\leq\b<+\infty$ and
$1<p<+\infty$ such that
\begin{equation}\label{pg1}
\b'|\xi|^p \leq W_\e(x,\xi) \leq \b(1+|\xi|^p)
\end{equation}
for a.e. $x \in \O$ and all $\xi \in \Rb^{3 \times 3}$. Define
$\J_\e : L^p(\O;\Rb^3) \times L^p(\o;\Rb^3) \times \A(\o) \to
[0,+\infty]$ by
$$\J_\e(u,\overline b,A):=\left\{
\begin{array}{ll}
\ds \int_{A \times I} W_\e\left(x,\nabla_\a
u\Big|\frac{1}{\e}\nabla_3 u\right)dx & \text{if }\left\{
\begin{array}{l}
u \in W^{1,p}(A \times I;\Rb^3),\\
\overline b = \frac{1}{\e}\int_I\nabla_3 u(\cdot,x_3)\, dx_3,
\end{array}\right.\\[0.3cm]
\ds +\infty & \text{otherwise.}
\end{array}
\right.$$

We prove the following integral representation for the $\G$-limit.

\begin{thm}\label{bending}
For every sequence $\{\e_n\} \searrow 0^+$, there exist a
subsequence (not relabeled) and a Carath\'eodory function $W^* : \o
\times \Rb^{3 \times 2} \times \Rb^3 \to [0,+\infty)$ (depending on
the subsequence) such that for every $A \in \A(\o)$, the sequence
$\J_{\e_n}(\cdot,\cdot,A)$ $\G$-converges in $L_s^p(A \times
I;\Rb^3) \times L_w^p(A;\Rb^3)$ to $\J(\cdot,\cdot,A)$ where
$$\J(u,\overline b,A)=\left\{
\begin{array}{ll}
\ds \int_A W^*(x_\a,\nabla_\a u(x_\a)|\overline b(x_\a))\, dx_\a &
\text{ if } u \in
W^{1,p}(A;\Rb^3),\\[0.3cm]
+\infty & \text{ otherwise}.
\end{array}
\right.$$
\end{thm}

\begin{proof}
For every $\{\e_n\} \searrow 0^+$, $u \in L^p(\O;\Rb^3)$, $\overline
b\in L^p(\o;\Rb^3)$ and $A \in \A(\o)$, let
$$\J(u,\overline b,A) :=  \inf_{\{u_n,\overline b_n\}}
\left\{\liminf_{n \to +\infty} \J_{\e_n}(u_n,\overline b_n,A) : u_n
\to u \text{ in }L^p(A \times I;\Rb^3) \text{ and } \overline b_n
\wto \overline b \text{ in } L^p(A;\Rb^3) \right\}.$$ Repeating word
for word the (standard) proof of \cite[Lemma~2.1]{BFMbis} one can
show that there exists a subsequence, still labeled $\{\e_n\}$, such
that for any $A \in \A(\o)$, $\J(\cdot,\cdot,A)$ is the $\G$-limit
in $L_s^p(A \times I;\Rb^3) \times L_w^p(A;\Rb^3)$ of
$\J_{\e_n}(\cdot,\cdot,A)$, that $\J(u,\overline b,A)=+\infty$ if $u
\in L^p(\O;\Rb^3) \setminus W^{1,p}(A;\Rb^3)$ and that for every
$(u,\overline b) \in W^{1,p}(\o;\Rb^3) \times L^p(\o;\Rb^3)$, the
set function $\J(u,\overline b,\cdot)$ is the restriction to
$\A(\o)$ of a Radon measure absolutely continuous with respect to
the Lebesgue measure $\LL^2$. The remaining of the proof is very
close to that of \cite[Theorem~1.1]{BDM}, thus we will only point
out the main changes. Let $\overline \xi \in \Rb^{3 \times 2}$, $z
\in \Rb^3$ and $x_0 \in \o$, define
$$W^*(x_0,\overline \xi|z):= \limsup_{\rho \to
0^+}\frac{\J(u_{\overline \xi},\overline
b_z,Q'(x_0,\rho))}{\rho^2}$$ where we have denoted $u_{\overline
\xi}(x_\a):=\overline\xi \,x_\a$ and $\overline b_z(x_\a):=z$. Since
$\J(u_{\overline \xi},\overline b_z,\cdot)$ is (the restriction of)
a Radon measure absolutely continuous with respect to $\LL^2$, we
have for every $A \in \A(\o)$,
\begin{equation}\label{1402}
\J(u_{\overline \xi},\overline b_z,A)= \int_A W^*(x_\a,\overline
\xi|z)\, dx_\a= \int_A W^*(x_\a,\nabla_\a u_{\overline
\xi}|\overline b_z)\, dx_\a.
\end{equation} By additivity, it is clear that
\begin{equation}\label{rep}
\J(u,\overline b,A)= \int_A W^*(x_\a,\nabla_\a u|\overline b)\,
dx_\a
\end{equation}
holds whenever $u$ is piecewise affine and $\overline b$ is
piecewise constant in $A$ and we wish to extend (\ref{rep}) to
arbitrary functions $u \in W^{1,p}(A;\Rb^3)$ and $\overline b \in
L^p(A;\Rb^3)$.

Using the lower semicontinuity of $\J$ and a suitable choice of
sequence, one can show as in \cite[Theorem~1.1]{BDM} that $\overline
\xi \mapsto W^*(x_0,\overline \xi|z)$ is rank one convex. We claim
that $z \mapsto W^*(x_0,\overline \xi|z)$ is convex. To see this let
$\theta \in [0,1]$, $z_1$, $z_2 \in \Rb^3$ and $\overline \xi \in
\Rb^{3 \times 2}$. Fix $x_0 \in \o$, $\rho>0$ and take an open set
$A \subset Q'(x_0,\rho)$ such that $\LL^2(\partial A)=0$ and
$\LL^2(A)=\theta \rho^2$ (take {\it e.g.}
$A=Q'(x_0,\sqrt{\theta}\rho)$). Define
$$\overline b_n(x_\a):=z_1 \chi(n x_\a) +z_2 (1-\chi(nx_\a))$$ where
$\chi$ is the characteristic function of $A$ in $Q'(x_0,\rho)$ which
has been extended to $\Rb^2$ by $\rho$-periodicity.
Riemann-Lebesgue's Lemma asserts that $\overline b_n \wto \overline
b_{\theta z_1 + (1-\theta)z_2}$ in $L^p(Q'(x_0,\rho);\Rb^3)$ and
since $\J(u_{\overline \xi},\cdot,Q'(x_0,\rho))$ is sequentially
weakly lower semicontinuous in $L^p(Q'(x_0,\rho);\Rb^3)$, it follows
that
\begin{eqnarray}\label{convexity}
\J(u_{\overline \xi},\overline b_{\theta z_1
+(1-\theta)z_2},Q'(x_0,\rho))& \leq & \liminf_{n \to
+\infty}\J(u_{\overline \xi},\overline
b_n,Q'(x_0,\rho))\nonumber\\
& = & \liminf_{n \to +\infty}\left\{\J(u_{\overline \xi},\overline
b_{z_1},A_n) + \J(u_{\overline \xi},\overline b_{z_2},Q'(x_0,\rho)
\setminus \overline A_n)\right\}
\end{eqnarray}
where $A_n:=\{x_\a \in Q'(x_0,\rho): \chi(nx_\a)=1\}$ is an open
set. Note that in the last equality, we have used the fact that
since $\LL^2(\partial A_n)=0$, then $\J(u_{\overline \xi},\overline
b_n,\partial A_n)=0$ as well and that $\J$ is local on open sets.
Using once more the Riemann-Lebesgue Lemma together with
(\ref{1402}), we get that
\begin{eqnarray*}\lim_{n \to +\infty}\J(u_{\overline \xi},\overline
b_{z_1},A_n) & = &  \lim_{n \to
+\infty}\int_{Q'(x_0,\rho)}\chi(nx_\a) W^*(x_\a,\overline
\xi|z_1)\, dx_\a\\
& = &  \theta \int_{Q'(x_0,\rho)} W^*(x_\a,\overline \xi|z_1)\,
dx_\a\\
& = & \theta \J(u_{\overline \xi},\overline b_{z_1},Q'(x_0,\rho))
\end{eqnarray*}
and similarly for the second term of (\ref{convexity}). Hence we
deduce that
$$\J(u_{\overline \xi},\overline
b_{\theta z_1 +(1-\theta)z_2},Q'(x_0,\rho)) \leq \theta
\J(u_{\overline \xi},\overline b_{z_1},Q'(x_0,\rho))
 + (1-\theta)\J(u_{\overline \xi},\overline b_{z_2},Q'(x_0,\rho))$$ and the convexity of
$W^*(x_0,\overline \xi|\cdot)$ arises after dividing the previous
inequality by $\rho^2$ and taking the $\limsup$ as $\rho$ tends to
zero. It follows that $(\overline \xi|z) \mapsto W^*(x_0,\overline
\xi|z)$ is separately convex for a.e. $x_0 \in \o$ and since the
following $p$-growth and $p$-coercivity conditions hold
\begin{equation}\label{pgW}
\b'(|\overline \xi|^p + |z|^p) \leq W^*(x_0,\overline \xi|z) \leq
\b(1+|\overline \xi|^p +|z|^p), \quad \text{ for a.e. }x_0 \in
\o\text{ and all }(\overline \xi,z) \in \Rb^{3 \times 2} \times
\Rb^3,
\end{equation}
we conclude that $(\overline \xi|z) \mapsto W^*(x_0,\overline
\xi|z)$ is continuous for a.e. $x_0 \in \o$ which proves that $W^*$
is a Carath\'eodory function.

We now prove that (\ref{rep}) holds for any $(u,\overline b)\in
W^{1,p}(A;\Rb^3) \times L^p(A;\Rb^3)$. By approximation and thanks
to the lower semicontinuity of $\J(\cdot,\cdot,A)$ for the strong
$W^{1,p}(A;\Rb^3) \times L^p(A;\Rb^3)$ topology, there holds
$$\J(u,\overline b,A) \leq \int_AW^*(x_\a,\nabla_\a u|\overline
b)\, dx_\a$$ for any $(u,\overline b)\in W^{1,p}(A;\Rb^3) \times
L^p(A;\Rb^3)$ and it remains to prove the converse inequality. This
is achieved exactly as in the final step of the proof of
\cite[Theorem~1.1]{BDM}, by considering the translated functional
$$\widetilde \J(v,\overline c,A):=\J(u+v,\overline b + \overline c,A)$$
where $(u,\overline b)$ are arbitrary functions in $W^{1,p}(A;\Rb^3)
\times L^p(A;\Rb^3)$.
\end{proof}

We refer to \cite{BF,BFMbis} for more explicit formulas for the
integrand $W^*$ in particular cases.\\

The following technical proposition states some kind of blow-up
result for functionals through $\G$-convergence. It will be of use
in the proof of the lower bound in Theorem {\ref{gammaconvbend}
because at some point, we will need to get rid of small residual
terms occurring inside the integrand $W_\e$. In \cite{BB1,BB2,babF},
this difficulty was treated thanks to a decoupling variable method
which consisted in replacing the function $W_\e$ by a much more
regular one thanks to Scorza-Dragoni's Theorem and Tietze's
Extension Theorem, and the set where these two integrands did not
match was controlled thanks to the equi-integrability result
\cite[Theorem~1.1]{BF}. This method was quite powerful in that
context since the manner on which $W_\e$ was depending on $\e$ was
completely known. However, in the generalized framework considered
here, it does not apply anymore since we have no information on the
way $W_\e$ depends on $\e$. The following blow up result, together
with a diagonalization argument (see Remark \ref{diag} below), will
enable us to overcome that problem.

\begin{proposition}\label{rho}
There exists a set $N \subset \o$ with $\LL^2(N)=0$ such that for
every $\{\rho_k\}\searrow 0^+$ and every $x_0 \in  \o \setminus N$,
the functional $J_k:L^p(B'\times I;\Rb^3) \times L^p(B';\Rb^3) \to
[0,+\infty]$ defined by
$$J_k(u,\overline b) =
\left\{\begin{array}{ll} \ds \int_{B'} W^*(x_0+\rho_k x_\a,\nabla_\a
u(x_\a)|\overline b(x_\a))\, dx_\a & \text{ if } u \in
W^{1,p}(B';\Rb^3)\\[0.3cm]
+\infty & \text{ otherwise},
\end{array}\right.$$
$\G$-converges in $L_s^p(B' \times I;\Rb^3) \times L_w^p(B';\Rb^3)$
to $J:L^p(B'\times I;\Rb^3) \times L^p(B';\Rb^3) \to [0,+\infty]$,
where
$$J(u,\overline b) =
\left\{\begin{array}{ll} \ds \int_{B'} W^*(x_0,\nabla_\a
u(x_\a)|\overline b(x_\a))\, dx_\a & \text{ if } u \in W^{1,p}(B';\Rb^3),\\[0.3cm]
+\infty & \text{ otherwise}.
\end{array}\right.$$
\end{proposition}

\begin{proof}
The proof relies on the Scorza-Dragoni Theorem (see {\it e.g.}
\cite[Chapter VIII]{ET}). For any $q \in \Nb$, there exists a
compact set $K_q \subset \o$ with $\LL^2(\o \setminus K_q) < 1/q$
and such that $W^*$ is continuous on $K_q \times \Rb^{3 \times 2}
\times \Rb^3$. Let $N:=\o \setminus \bigcup_q K_q^*$ where
\begin{equation}\label{measurebis} K_q^*:=\left\{x \in K_q :
\lim_{\rho \to 0} \frac{\LL^2(B'(x_0,\rho) \setminus
K_q)}{\LL^2(B'(x_0,\rho))} = 0\right\}.
\end{equation}
Since $\LL^2(K_q \setminus K_q^*)=0$, then $\LL^2(N) \leq \LL^2 (\o
\setminus K^*_q)=\LL^2(\o \setminus K_q)<1/q \to 0$. Select a point
$x_0 \in \o \setminus N$, so that $x_0 \in K_q^*$ for some $q \in
\Nb$.

{\bf The upper bound. }Assume first that $u \in
W^{1,\infty}(B';\Rb^3)$ and $\overline b \in L^\infty(B';\Rb^3)$ and
set $M:= \|(\nabla_\a u | \overline b) \|_{L^\infty(B';\Rb^{3 \times
3})}$. Then according to the $p$-growth condition (\ref{pgW})
\begin{eqnarray}\label{wiinf}
J_k(u,\overline b) & = & \int_{B'} W^*(x_0 + \rho_k x_\a ,\nabla_\a u|\overline b)\, dx_\a\nonumber\\
& \leq & \int_{B' \cap \left( \frac{K_q - x_0}{\rho_k}\right)}
W^*(x_0 + \rho_k x_\a ,\nabla_\a u|\overline b)\, dx_\a +\b (1+M^p)
\LL^2 \left(B' \setminus \left(\frac{K_q -
x_0}{\rho_k}\right)\right).
\end{eqnarray}
As $W^*$ is uniformly continuous on $K_q \times B(0,M)$, there
exists a continuous and increasing function $\eta :[0,+\infty) \to
[0,+\infty)$ such that $\eta(0)=0$ and
\begin{equation}\label{unifcont}
\int_{B' \cap \left( \frac{K_q - x_0}{\rho_k}\right)}|W^*(x_0 +
\rho_k x_\a ,\nabla_\a u|\overline b) - W^*(x_0, \nabla_\a
u|\overline b)|\, dx_\a \leq \eta(\rho_k).
\end{equation}
Gathering (\ref{measurebis}), (\ref{wiinf}) and (\ref{unifcont}) and
passing to the limit as $k \to +\infty$ yields to
$$\G\text{-}\limsup_{k \to +\infty}J_k(u,\overline b) \leq
\limsup_{k \to +\infty} J_k(u,\overline b) \leq  J(u,\overline b).$$
The general case follows from the density of $W^{1,\infty}(B';\Rb^3)
\times L^\infty(B';\Rb^3)$ in $W^{1,p}(B';\Rb^3) \times
L^p(B';\Rb^3)$, the lower continuity of the $\G$-limsup and the
continuity of $J$ for the strong $W^{1,p}(B';\Rb^3) \times
L^p(B';\Rb^3)$-topology.

{\bf The lower bound. }Let $(u,\overline b) \in L^p(B' \times
I;\Rb^3) \times L^p(B';\Rb^3)$ and $\{(u_k,\overline b_k)\} \subset
L^p(B' \times I;\Rb^3) \times L^p(B';\Rb^3)$ such that $u_k \to u$
in $L^p(B'\times I;\Rb^3)$, $\overline b_k \wto \overline b$ in
$L^p(B';\Rb^3)$ and
$$\liminf_{k \to +\infty}J_k(u_k,\overline b_k)<+\infty.$$
Up to a subsequence (not relabeled) we can suppose that $u$ and $u_k
\in W^{1,p}(B';\Rb^3)$ for each $k\in \Nb$ and that $u_k \wto u$ in
$W^{1,p}(B';\Rb^3)$. According to the Decomposition Lemma
\cite[Lemma~1.2]{FMP} and Chacon's Biting Lemma
\cite[Lemma~5.32]{AFP}, there is no loss of generality to assume
that $\{|\nabla_\a u_k|^p\}$ and $\{|\overline b_k|^p\}$ are
equi-integrable. Define the set $A_k^t:=\left\{x_\a \in B':
|(\nabla_\a u_k(x_\a)|\overline b_k(x_\a))|\leq t\right\}$. From
Chebyshev's Inequality we have that $\LL^2(B'\setminus A_k^t) \leq c
/ t^p$ for some constant $c>0$ independent of $t$ and $k$ and
arguing exactly as in the proof of the upper bound, one can show
that for each $t>0$,
\begin{eqnarray}\label{wsup}
\liminf_{k \to +\infty}J_k(u_k,\overline b_k) \geq \liminf_{k \to
+\infty}\int_{A_k^t \cap \left( \frac{K_q - x_0}{\rho_k}\right)}
W^*(x_0,\nabla_\a u_k|\overline b_k)\, dx_\a.
\end{eqnarray}
According to the $p$-growth condition (\ref{pgW}) and
(\ref{measurebis}),
\begin{equation}\label{Kj}
\int_{A_k^t \setminus \left(\frac{K_q - x_0}{\rho_k}\right)}
W^*(x_0,\nabla_\a u_k|\overline b_k)\, dx_\a \leq \b(1+t^p) \LL^2
\left( B' \setminus \left( \frac{K_q - x_0}{\rho_k}\right)\right)
\xrightarrow[k \to +\infty]{} 0
\end{equation}
while the equi-integrability of $\{|\nabla_\a u_k|^p\}$ and
$\{|\overline b_k|^p\}$ and the fact that $\LL^2(B' \setminus A_k^t)
\to 0$ as $t \to +\infty$ imply that
\begin{equation}\label{AkM}
\sup_{k \in \Nb}\int_{B' \setminus A_k^t} W^*(x_0,\nabla_\a
u_k|\overline b_k)\, dx_\a \leq \b \sup_{k \in \Nb} \int_{B'
\setminus A_k^t}(1+|\nabla_\a u_k|^p+|\overline b_k)|^p)\, dx_\a
\xrightarrow[t \to +\infty]{} 0.
\end{equation}
Hence gathering (\ref{wsup}), (\ref{Kj}) and (\ref{AkM}) yields to
$$\liminf_{k \to +\infty} J_k(u_k,\overline b_k) \geq
\liminf_{k \to +\infty} J(u_k,\overline b_k) \geq J(u,\overline b)$$
where the last inequality holds because $J$ is sequentially weakly
lower semicontinuous in $W^{1,p}(B';\Rb^3) \times L^p(B';\Rb^3)$.
\end{proof}

\begin{rmk}\label{diag}{\rm
One can show that in Theorem \ref{bending}, the value of $\J$ does
not change replacing $W_{\e_n}$ by is quasiconvexification
$QW_{\e_n}$ defined by
\begin{equation}\label{quasiconvexification}
QW_{\e_n}(x,\xi):= \inf_{\varphi \in W^{1,\infty}_0((0,1)^3;\Rb^3)}
\int_{(0,1)^3} W_{\e_n}(x,\xi + \nabla \varphi(y))\, dy \quad \text{
for all }\xi \in \Rb^{3 \times 3}\text{ and a.e. }x\in \O.
\end{equation}
Hence there is no loss of generality to assume in Theorem
\ref{bending} that $W_\e$ is quasiconvex. Since the weak topology on
every normed bounded subsets of $L^p(B';\Rb^3)$ is metrizable, it
follows from a diagonalization argument, Theorem \ref{bending},
Proposition \ref{rho} and the fact that $\G$-convergence of coercive
and lower semicontinuous functionals on a metric space is metrizable
(see \cite[Theorem~10.22~(a)]{DM}), that for every $M>0$ and every
sequence $\{\rho_k\} \searrow 0^+$, there exists a subsequence $n(k)
\nearrow +\infty$ such that $\e_{n(k)}/\rho_k \to 0$ and for every
$(u,\overline b) \in L^p(B'\times I;\Rb^3)\times L^p(B';\Rb^3)$ with
$\|\overline b\|_{L^p(B';\Rb^3)} \leq M$, then the $\G$-limit in
$L_s^p(B' \times I;\Rb^3) \times L_w^p(B';\Rb^3)$ of
$$\left\{\begin{array}{ll} \ds \int_{B'\times I}
W_{\e_{n(k)}}\left(x_0+\rho_k x_\a,x_3,\nabla_\a
u\Big|\frac{\rho_k}{\e_{n(k)}}\nabla_3 u\right) dx & \text{ if }
\left\{
\begin{array}{l}
u \in W^{1,p}(B'\times I;\Rb^3),\\
\overline b=\frac{\rho_k}{\e_{n(k)}}\int_I \nabla_3 u(\cdot,x_3)\,
dx_3,
\end{array}
\right.\\[0.3cm]
+\infty & \text{ otherwise},
\end{array}\right.$$
coincides with
$$\left\{\begin{array}{ll} \ds \int_{B'} W^*(x_0,\nabla_\a
u|\overline
b)\, dx_\a & \text{ if } u \in W^{1,p}(B';\Rb^3),\\[0.3cm]
+\infty & \text{ otherwise}
\end{array}\right.$$
for every $x_0 \in \o \setminus N$, where $N \subset \o$ is the same
exceptional set than in Proposition \ref{rho}. }
\end{rmk}

\section{Integral representation for dimension reduction problems in
$SBV$ involving the bending moment}\label{gammaconvergence}

\noindent We now come to the heart of this study that is dealing
with a similar problem than in Theorem \ref{bending} but in the
framework of Special functions with Bounded Variation, adding a
surface energy term. Let us define $\mathcal G_\e : BV(\O;\Rb^3)
\times L^p(\o;\Rb^3) \to [0,+\infty]$ by
$$\mathcal G_\e(u,\overline b):=\left\{
\begin{array}{ll}
\ds \int_\O W_\e\left(x,\nabla_\a u\Big|\frac{1}{\e}\nabla_3
u\right)dx + \int_{S_u} \left| \left( \left(\nu_u\right)_\a
\Big|\frac{1}{\e} \left( \nu_u \right)_3 \right)\right| d\mathcal
H^2 & \text{if } \left\{ \begin{array}{l}
u \in SBV^p(\O;\Rb^3),\\
\overline b = \frac{1}{\e}\int_I\nabla_3 u(\cdot,x_3)\, dx_3,
\end{array} \right.\\[0.3cm]
\ds +\infty & \text{otherwise.}
\end{array}
\right.$$

\noindent Then, the following $\G$-convergence result holds:

\begin{thm}\label{gammaconvbend}
For every sequence $\{\e_n\} \searrow 0^+$, there exists a
subsequence, still labeled $\{\e_n\}$ such that $\mathcal G_{\e_n}$
$\G$-converges in $L_s^1(\O;\Rb^3) \times L_w^p(\o;\Rb^3)$ to
$\mathcal G : BV(\O;\Rb^3) \times L^p(\o;\Rb^3) \to [0,+\infty]$
defined by
$$\mathcal G(u,\overline b):=\left\{
\begin{array}{ll}
\ds \int_\o W^*(x_\a,\nabla_\a u|\overline b)\, dx_\a + \HH^1(S_u) & \text{if }  u \in SBV^p(\o;\Rb^3),\\[0.3cm]
\ds  +\infty & \text{otherwise},
\end{array}
\right.$$ where $W^*$ is given by Theorem \ref{bending}.
\end{thm}

The remaining of this section is devoted to prove Theorem
\ref{gammaconvbend}. We will first localize the functional $\mathcal
G_\e$ on $\A(\o)$, and noticing that minimizing sequences are not
necessarily weakly relatively compact in $BV$, we will use the same
truncation argument than in \cite{babadjian} (see also \cite{FF})
introducing an artificial functional. Then we will show that it
actually coincides with the $\G$-limit whenever $u \in BV(\O;\Rb^3)
\cap L^\infty(\O;\Rb^3)$ (see Lemma \ref{bd} and Remark \ref{infty})
and it will enable us to show that for such $u$'s the $\G$-limit is
a measure absolutely continuous with respect to $\LL^2+\HH^1 \res\,
S_u$ (see Lemma \ref{measure}). Together with a blow up argument,
this property will be useful to prove the upper bound in Lemma
\ref{upper} while the lower bound, Lemma \ref{lower}, will obtained
thanks to Theorem \ref{bab} and a suitable diagonalization argument
(see Remark \ref{diag}).

\subsection{Localization}

We first localize our functional on $\A(\o)$ defining $\mathcal G_\e
: BV(\O;\Rb^3) \times L^p(\o;\Rb^3) \times \A(\o) \to [0,+\infty]$
by
$$\mathcal G_\e(u,\overline b,A):=\left\{
\begin{array}{ll}
\begin{array}{l} \ds \int_{A \times I}  \ds W_\e\left(x,\nabla_\a u\Big|\frac{1}{\e}\nabla_3 u\right)dx  \\
\hspace{1cm} \ds + \int_{S_u \cap (A \times I)} \left| \left(
\left(\nu_u\right)_\a \Big|\frac{1}{\e} \left( \nu_u \right)_3
\right)\right| d\mathcal H^2 \end{array} & \text{if }\left\{
\begin{array}{l}
u \in SBV^p(A \times I;\Rb^3),\\
\overline b = \frac{1}{\e}\int_I\nabla_3 u(\cdot,x_3)\, dx_3,
\end{array}\right.\\
\ds +\infty & \text{ otherwise}.
\end{array}
\right.$$ For every sequence $\{\e_n\} \searrow 0^+$ and all
$(u,\overline b,A) \in BV(\O;\Rb^3) \times L^p(\o;\Rb^3) \times
\A(\o)$, we define
\begin{equation}\label{IA}
\E (u,\overline b,A):= \inf_{\{u_n,\overline b_n\}} \left\{
\liminf_{n \to +\infty} \; \mathcal G_{\e_n}(u_n,\overline b_n,A) :
u_n \to u \text{ in }L^1(A \times I;\Rb^3), \, \overline b_n \wto
\overline b \text{ in } L^p(A;\Rb^3) \right\}.
\end{equation}
Theorem 8.5 and Corollary 8.12 in \cite{DM} together with a
diagonalization argument imply  the existence of a subsequence,
still denoted $\{\e_n\},$ such that, for any $A \in \mathcal R(\o)$
(or $A=\o$), $\E(\cdot,\cdot,A)$ is the $\G$-limit of $\mathcal
G_{\e_n}(\cdot,\cdot,A)$ in $L_s^1(A \times I;\Rb^3) \times
L_w^p(A;\Rb^3)$. Extracting if necessary a further subsequence, one
may assume that $\{\e_n\}$ is chosen so that Theorem \ref{bending}
holds. To prove Theorem \ref{gammaconvbend}, it is enough to show
that $\E(u,\overline b,\o)=\mathcal G(u,\overline b)$.


\subsection{A truncation argument}\label{truncation}

As pointed out in \cite{babadjian}, the main problem with the
definition of $\E$ in (\ref{IA}) is that minimizing sequences are
not necessarily bounded in $BV(\O;\Rb^3)$ and thus, not necessarily
weakly convergent in this space. Following \cite{babadjian}, we
define for all $(u,\overline b,A) \in BV(\O;\Rb^3) \times
L^p(\o;\Rb^3) \times \A(\o)$
\begin{eqnarray*}\E_\infty (u,\overline b,A)&:=&
\inf_{\{u_n,\overline b_n\}} \Big\{ \liminf_{n \to +\infty} \;
\mathcal G_{\e_n}(u_n,\overline b_n,A) : u_n \to u \text{ in
}L^1(A \times I;\Rb^3),\\
&&\hspace{1.2cm} \overline b_n \wto \overline b \text{ in
}L^p(A;\Rb^3),\; \sup_{n \in \Nb} \|u_n\|_{L^\infty(A \times
I;\Rb^3)} <+\infty \Big\}. \end{eqnarray*}It is immediate that
$\E(u,\overline b,A) \leq \E_\infty (u,\overline b,A)$ while we will
show that equality holds when $u$ belongs to $BV(\O;\Rb^3) \cap
L^\infty(\O;\Rb^3)$. This will be obtained as a consequence of Lemma
\ref{bd} below. It means that for such deformation fields $u \in
BV(\O;\Rb^3) \cap L^\infty(\O;\Rb^3)$, strong
$L^1(\O;\Rb^3)$-convergence and weak $BV(\O;\Rb^3)$-convergence are,
in a sense, equivalent for the computation of the $\G$-limit.

\begin{lemma}\label{bd}
Let $A \in \A(\o)$, $u \in BV(\O;\Rb^3)\cap L^\infty(\O;\Rb^3)$ and
$\overline b \in L^p(\o;\Rb^3)$. If $\{u_n\} \subset SBV^p(A \times
I;\Rb^3)$ is such that $u_n \to u$ in $L^1(A \times I;\Rb^3)$,
$\frac{1}{\e_n}\int_I \nabla_3 u_n(\cdot,x_3)\, dx_3 \wto \overline
b$ in $L^p(A;\Rb^3)$ and the following limit
$$L:=\lim_{n \to +\infty}  \mathcal
G_{\e_n}\left(u_n,\frac{1}{\e_n}\int_I\nabla_3 u_n(\cdot,x_3)\,
dx_3,A\right)$$ exists and is finite. Then, for any $\eta>0$ one can
find $C>0$ and $\{w_n\} \subset SBV^p(A \times I;\Rb^3)$ such that
$w_n \to u$ in $L^1(A \times I;\Rb^3)$,
$\frac{1}{\e_n}\int_I\nabla_3 w_n(\cdot,x_3)\, dx_3 \wto \overline
b$ in $L^p(A;\Rb^3)$, $\sup_n \|w_n\|_{L^\infty(A \times I;\Rb^3)}
\leq C$ and
$$L \geq \limsup_{n \to +\infty} \mathcal
G_{\e_n}\left(w_n,\frac{1}{\e_n}\int_I\nabla_3 w_n(\cdot,x_3)\,
dx_3,A\right)-\eta.$$
\end{lemma}

\begin{proof} Let us define a smooth truncation function $\varphi_i \in
\mathcal C^1_c(\Rb^3 ; \Rb^3)$ satisfying
\begin{equation}\label{fii}
\varphi_i(s)=\left\{
\begin{array}{rcl}
s & \text{if} & |s|<e^i,\\[0.2cm]
0 & \text{if} & |s| \geq e^{i+1}
\end{array}\right. \quad \text{ and }\quad |\nabla \varphi_i(s)|\leq 2.
\end{equation}
Let $w_{n,i}:=\varphi_i(u_{n})$, thanks to the Chain Rule formula
\cite[Theorem~3.96]{AFP}, $w_{n,i} \in SBV^p(A \times I;\Rb^3)$ and
\begin{equation}\label{wik}
\left\{\begin{array}{l}
\|w_{n,i}\|_{L^\infty(A \times I;\Rb^3)} \leq e^{i+1},\\[0.2cm]
S_{w_{n,i}} \subset S_{u_n},\\[0.2cm]
\nabla w_{n,i} = \nabla \varphi_i(u_n) \nabla u_{n} \quad \mathcal
L^3\text{-a.e. in }A \times I.
\end{array}\right.
\end{equation}
Since $u \in L^\infty(\O;\Rb^3)$, we can choose $i$ large enough ($i
\geq m:=[\ln(\|u\|_{L^\infty(\O;\Rb^3)})]+1$) so that
$u=\varphi_i(u)$ and thus according to (\ref{fii})
\begin{equation}\label{L1}
\|w_{n,i}-u\|_{L^1(A \times I;\Rb^3)}  =
\|\varphi_i(u_{n})-\varphi_i(u)\|_{L^1(A \times I;\Rb^3)} \leq 2
\|u_{n}-u\|_{L^1(A \times I;\Rb^3)}.
\end{equation}
Since (a subsequence of) $u_n \to u$ a.e. in $A \times I$ and
$\nabla \varphi_i$ is continuous, it follows that $\nabla
\varphi_i(u_n) \to \nabla \varphi_i(u)={\rm Id}$ a.e. in $A\times I$
as $n \to +\infty$. Take $v \in L^{p'}(A;\Rb^3)$ where $1/p+1/p'=1$,
as $|\nabla \varphi_i(u_n)^T v| \leq 2|v| \in L^{p'}(A)$, the
Dominated Convergence Theorem implies that $\nabla \varphi_i(u_n)^T
v \to v$ in $L^{p'}(A \times I;\Rb^3)$ and thus
\begin{eqnarray*}\lim_{n \to +\infty}\int_A \left(\frac{1}{\e_n}\int_I \nabla_3
w_{n,i}(x_\a,x_3)\,dx_3 \right) \cdot v(x_\a)\, dx_\a & = & \lim_{n
\to +\infty}\int_{A \times I}\frac{1}{\e_n}\nabla_3 u_n \cdot\left(
\nabla \varphi_i(u_n)^T v \right)\,
dx\\
& = & \int_A \overline b \cdot v\, dx_\a,
\end{eqnarray*}
where we used the fact that $(1/\e_n)\nabla_3 u_n \rightharpoonup b$
in $L^p(A \times I;\Rb^3)$ and $\overline b=\int_I b(\cdot,x_3)\,
dx_3$. Hence
\begin{equation}\label{weakconv}
\frac{1}{\e_n}\int_I \nabla_3 w_{n,i}(\cdot,x_3)\, dx_3
\xrightharpoonup[n \to +\infty]{} \overline b \text{ in
}L^p(A;\Rb^3), \text{ for all }i \geq m.\end{equation} The growth
condition (\ref{pg1}), (\ref{fii}) and (\ref{wik}) imply that
\begin{eqnarray}\label{sous}
&&\int_{A \times I}W_{\e_n} \left(x, \nabla_\a w_{n,i} \Big|\frac{1}{\e_{n}} \nabla_3 w_{n,i}\right)dx\nonumber\\
&&\hspace{2.0cm}\leq \int_{\{|u_n| < e^i\} }W_{\e_n}\left(x,
\nabla_\a u_{n} \Big|\frac{1}{\e_{n}} \nabla_3 u_{n} \right)dx
+ \b \mathcal L^3(\{|u_n| \geq e^{i+1}\})\nonumber\\
&&\hspace{2.5cm} +\int_{\{e^i \leq  |u_n| < e^{i+1}\} }W_{\e_n}
\left(x, \nabla \varphi_i(u_n) \nabla_\a u_{n} \Big|
\frac{1}{\e_{n}} \nabla \varphi_i(u_n) \nabla_3 u_{n} \right)dx\nonumber\\
&&\hspace{2.0cm} \leq \int_{A \times I}W_{\e_n}
\left(x, \nabla_\a u_{n} \Big|\frac{1}{\e_{n}} \nabla_3 u_{n} \right)dx +
\beta\, e^{-i}\, \|u_{n}\|_{L^1(A \times I;\Rb^3)}\nonumber\\
&&\hspace{2.5cm} + 2^p \beta  \int_{\{e^i \leq |u_n| < e^{i+1}\}
}\left|\left( \nabla_\a u_{n} \Big|\frac{1}{\e_{n}}  \nabla_3 u_{n}
\right)\right|^p dx,
\end{eqnarray}
where we have used Chebyshev's Inequality. Since
$\nu_{w_{n,i}}(x)=\pm \nu_{u_{n}}(x)$ for $\mathcal H^2$-a.e. $x \in
S_{w_{n,i}}$, (\ref{wik}) yields to
\begin{equation}\label{sur}
\int_{S_{w_{n,i}} \cap (A \times I) } \left| \left(
\left(\nu_{w_{n,i}}\right)_\a \Big| \frac{1}{\e_{n}}
\left(\nu_{w_{n,i}}\right)_3 \right)\right| d\mathcal H^2 \leq
\int_{S_{u_{n}} \cap (A \times I)} \left| \left(
\left(\nu_{u_{n}}\right)_\a \Big| \frac{1}{\e_{n}}
\left(\nu_{u_{n}}\right)_3 \right)\right| d\mathcal H^2.
\end{equation}
Let $M \in \Nb$, from (\ref{sous}) and (\ref{sur}), a summation for
$i=m$ to $M$ implies that
\begin{eqnarray*}
&&\frac{1}{M-m+1} \sum_{i=m}^M \left[ \int_{A \times I}
W_{\e_n}\left(x, \nabla_\a w_{n,i} \Big|\frac{1}{\e_{n}} \nabla_3 w_{n,i}\right)dx\right.\\
&&\hspace{6cm}\left.+ \int_{S_{w_{n,i}} \cap (A \times I) } \left|
\left( \left(\nu_{w_{n,i}}\right)_\a
\Big| \frac{1}{\e_{n}} \left(\nu_{w_{n,i}}\right)_3 \right)\right| d\mathcal H^2 \right]\\
&&\hspace{0.5cm} \leq \int_{A \times I} W_{\e_n}\left(x, \nabla_\a
u_{n} \Big|\frac{1}{\e_{n}} \nabla_3 u_{n} \right)dx +
\int_{S_{u_{n}} \cap (A \times I)} \left| \left(
\left(\nu_{u_{n}}\right)_\a \Big| \frac{1}{\e_{n}}
\left(\nu_{u_{n}}\right)_3 \right)\right| d\mathcal H^2+
\frac{c}{M-m+1},
\end{eqnarray*}
where $$c=\beta\sup_{n \in \Nb}\|u_{n}\|_{L^1(A \times I;\Rb^3)}
\sum_{i \geq 1}e^{-i} + 2^p \beta \sup_{n \in \Nb}\left\|\left(
\nabla_\a u_{n} \Big|\frac{1}{\e_{n}}  \nabla_3 u_{n} \right)
\right\|^p_{L^p(A \times I;\Rb^{3 \times 3})} < +\infty.$$ We may
find some $i_n \in \{m,\ldots,M\}$ such that, setting
$w_n:=w_{n,i_n}$, then
\begin{eqnarray}\label{last}
&&\hspace{-1cm}\int_{A \times I}W_{\e_n}\left(x,\nabla_\a w_n
\Big|\frac{1}{\e_{n}} \nabla_3 w_n\right)dx + \int_{S_{w_n} \cap (A
\times I)} \left| \left( \left(\nu_{w_n}\right)_\a
\Big| \frac{1}{\e_{n}} \left(\nu_{w_n}\right)_3 \right)\right| d\mathcal H^2 \nonumber\\
&& \hspace{0cm}\leq \int_{A \times I}W_{\e_n}\left(x,\nabla_\a u_{n}
\Big|\frac{1}{\e_{n}} \nabla_3 u_{n} \right)dx + \int_{S_{u_{n}}
\cap (A \times I)} \left| \left( \left(\nu_{u_{n}}\right)_\a \Big|
\frac{1}{\e_{n}} \left(\nu_{u_{n}}\right)_3 \right)\right| d\mathcal
H^2+ \frac{c}{M-m+1}.
\end{eqnarray}
Moreover, in view of (\ref{L1}) and (\ref{weakconv}), $w_n \to u$ in
$L^1(A \times I;\Rb^3)$, $\frac{1}{\e_n}\int_I \nabla_3
w_n(\cdot,x_3)\, dx_3 \wto \overline b$ in $L^p(A;\Rb^3)$ and
(\ref{wik}) implies that $\|w_n\|_{L^\infty(A \times I;\Rb^3)} \leq
e^{i_n+1} \leq e^{M+1}$. The proof is achieved passing to the limit
as $n$ tends to $+\infty $ in (\ref{last}) and choosing $M$ large
enough so that $c/(M-m+1) \leq \eta$.
\end{proof}

\begin{rmk}\label{infty}{\rm
As a consequence of Lemma \ref{bd}, we get that for any $A \in
\mathcal R(\o)$ (or $A=\o$), every $u \in BV(\O;\Rb^3) \cap
L^\infty(\O;\Rb^3)$ and every $\overline b \in L^p(\o;\Rb^3)$, then
$\E(u,\overline b,A)=\E_\infty(u,\overline b,A)$. }\end{rmk}

\begin{rmk}\label{sobo}{\rm
A similar statement of Lemma \ref{bd} can be proved in the framework
of Sobolev spaces, replacing $\mathcal G_{\e_n}$ by $\J_{\e_n}$.
}\end{rmk}

\begin{rmk}\label{WQW}
{\rm Using a relaxation argument in $SBV^p$ as in the proof of
\cite[Lemma~3.4]{babadjian} and Lemma \ref{bd}, one can show that if
$u \in SBV^p(\o;\Rb^3)\cap L^\infty(\o;\Rb^3)$ and if $\overline b
\in L^p(\o;\Rb^3)$, the value of $\E_\infty$ does not change
replacing $W_{\e_n}$ by is quasiconvexification $QW_{\e_n}$ defined
in (\ref{quasiconvexification}). The main point is that the
diagonalization argument can still be used despite the weak
$L^p(\o;\Rb^3)$-convergence of the bending moment since the dual of
$L^p(\o;\Rb^3)$ is separable. Hence we may assume without loss of
generality that $W_\e$ is quasiconvex. In particular (see
\cite[Lemma~2.2,~Chapter~4]{D}), the following $p$-Lipschitz
condition holds,
\begin{equation}\label{pg1bis}
|W_{\e}(x,\xi_1)- W_{\e}(x,\xi_2)| \leq c
(1+|\xi_1|^{p-1}+|\xi_2|^{p-1})|\xi_1-\xi_2|, \text{ for all }
\xi_1,\; \xi_2 \in \Rb^{3 \times 3}\text{ and a.e. }x\in \O.
\end{equation}
}
\end{rmk}


Lemma \ref{bd} and Remark \ref{infty} are essential for the proof of
the following result because they allow us to replace strong
$L^1(\O;\Rb^3)$-convergence of any minimizing sequence by strong
$L^p(\O;\Rb^3)$-convergence.

\begin{lemma}\label{measure}
For all $u \in SBV^p(\o;\Rb^3) \cap L^\infty(\o;\Rb^3)$ and all
$\overline b \in L^p(\o;\Rb^3)$, $\E_\infty(u,\overline b,\cdot)$ is
the restriction to $\A(\o)$ of a Radon measure absolutely continuous
with respect to $\mathcal L^2 + \mathcal H^1\res \,  S_u$.
\end{lemma}

\begin{proof}
Let $u \in SBV^p(\o;\Rb^3) \cap L^\infty(\o;\Rb^3)$, $A \in \A(\o)$
and assume first that $\overline b$ is smooth. Then taking
$u_n(x_\a,x_3):=u(x_\a)+\e_n x_3 \overline b(x_\a)$ and $\overline
b_n(x_\a):=\overline b(x_\a)$ as test functions for
$\E_\infty(u,\overline b,A)$ and using the $p$-growth condition
(\ref{pg1}), we get that
\begin{equation}\label{harraps}
\E_\infty(u,\overline b,A) \leq \b \int_{A}(1+|\nabla_\a u|^p +
|\overline b|^p)\, dx_\a+ \HH^1(S_u\cap A). \end{equation} The same
inequality holds for arbitrary functions $\overline b \in
L^p(\o;\Rb^3)$ thanks to the density of smooth maps into
$L^p(\o;\Rb^3)$ and the sequential weak lower semicontinuity of
$\E_\infty(u,\cdot,A)$ in $L^p(A;\Rb^3)$. The remaining of the proof
is very classical and is essentially the same than that of
\cite[Lemma~3.6]{babadjian}. As usual, the most delicate point is to
prove the subadditivity of $\E_\infty(u,\overline b,\cdot)$ and this
is done by gluing together suitable minimizing sequences by means of
a cut-off function. The argument still works with the presence of
the bending moment since the cut-off function is chosen
independently of $x_3$. One should once more be careful when
applying a diagonalization argument because of the weak convergence
in $L^p$. As already mentioned in Remark \ref{WQW}, it is still
allowed in the case where we include the bending moment since dual
of $L^p$ is separable.
\end{proof}

As a consequence of Lemma \ref{measure} and Lebesgue's Decomposition
Theorem, there exists a $\mathcal L^2$-measurable function $f$ and a
$\mathcal H^1\res \,  S_u$-measurable function $g$ such that for
every $A \in \A(\o)$,
\begin{equation}\label{intrep}
\E_\infty(u,\overline b,A)= \int_A f\, d\mathcal L^2 + \int_{A \cap
S_u} g\, d\mathcal H^1.
\end{equation}Since the measures
$\mathcal L^2$ and $\mathcal H^1\res \,  S_u$ are mutually singular,
$f$ is the Radon-Nikod\'ym derivative of $\E_\infty(u,\overline
b,\cdot)$ with respect to $\mathcal L^2$,
$$f(x_0)=\lim_{\rho \to 0}\frac{\E_\infty(u,\overline b,B'(x_0,\rho))}{\LL^2(B'(x_0,\rho))},
\quad \text{for } \mathcal L^2 \text{-a.e. } x_0 \in \o$$ and $g$ is
the Radon-Nikod\'ym derivative of $\E_\infty(u,\overline b,\cdot)$
with respect to $\mathcal H^1\res \,  S_u$,
$$g(x_0)=\lim_{\rho \to 0}\frac{\E_\infty(u,\overline b,B'(x_0,\rho))}{\mathcal H^1(S_u \cap B'(x_0,\rho))},
\quad \text{for } \mathcal H^1 \text{-a.e. } x_0 \in S_u.$$


\subsection{The upper bound}\label{upper}

\noindent We first show the upper bound. To this end, we will use
the locality property of the $\G$-limit proved in the previous
subsection when $u \in BV(\O;\Rb^3) \cap L^\infty(\O;\Rb^3)$ and the
analogue $\G$-convergence result in Sobolev spaces (Theorem
\ref{bending}).

\begin{lemma}\label{gs1}
For all $u \in BV(\O;\Rb^3)$ and all $\overline b \in
L^p(\o;\Rb^3)$, $\E(u,\overline b,\o) \leq \mathcal G(u,\overline
b)$.
\end{lemma}

\begin{proof}It is enough to consider the case where $\mathcal G(u,\overline
b)<+\infty$ and thus $u \in SBV^p(\o;\Rb^3)$. In fact, we will first
restrict to the case where $u \in L^\infty(\o;\Rb^3) \cap
SBV^p(\o;\Rb^3)$ because thanks to Remark \ref{infty}, it allows us
to replace $\E$ by $\E_\infty$. According to (\ref{intrep}) and the
definition of $\mathcal G$, we must show that $g(x_0) \leq 1$ for
$\mathcal H^1$-a.e. $x_0 \in S_u$ and
$f(x_0) \leq W^* (x_0,\nabla_\a u(x_0)|\overline b(x_0))$ for $\mathcal L^2$-a.e. $x_0 \in \o$.\\

Let us first treat the surface term. By virtue of (\ref{harraps})
with $A=B'(x_0,\rho)$, we have that for $\HH^1$-a.e. $x_0 \in S_u$,
\begin{eqnarray*}
g(x_0) & = & \lim_{\rho \to 0}\frac{\E_\infty(u,\overline b,B'(x_0,\rho))}{\mathcal H^1(S_u \cap B'(x_0,\rho))}\\
& \leq & \limsup_{\rho \to 0} \frac{1}{\mathcal H^1(S_u \cap
B'(x_0,\rho))}\left\{\b \int_{B'(x_0,\rho)}
(1+|\nabla_\a u|^p+|\overline b|^p) \,dx_\a +  \mathcal H^1(S_u \cap B'(x_0,\rho)) \right\}\\
& =  & \limsup_{\rho \to 0} \frac{\mu(B'(x_0,\rho))}{\mathcal
H^1(S_u \cap B'(x_0,\rho))} + 1,
\end{eqnarray*}
where we set $\mu:= \b(1+|\nabla_\a u|^p+|\overline b|^p) \mathcal
L^2$. But since $\mu$ and $\mathcal H^1\res \,  S_u$ are mutually
singular, we have for $\mathcal H^1$-a.e. $x_0 \in S_u$
$$\lim_{\rho \to 0}\frac{\mu(B'(x_0,\rho))}{\mathcal H^1(S_u \cap B'(x_0,\rho))} = 0,$$
which shows that $g(x_0) \leq 1$ for $\mathcal H^1$-a.e. $x_0 \in S_u$.\\

Concerning the bulk term, choose $x_0 \in \o$ to be a Lebesgue point
of $u$, $\nabla_\a u$, $\overline b$ and $W^*(\cdot,\nabla_\a
u(\cdot)|\overline b(\cdot))$ and such that
\begin{equation}\label{sing}
\lim_{\rho \to 0} \frac{\mathcal H^1(S_u \cap
B'(x_0,\rho))}{\LL^2(B'(x_0,\rho))}=0.
\end{equation}Remark that $\mathcal
L^2$ almost every points $x_0$ in $\omega$ satisfy these properties
and set $u_0(x_\a):=\nabla_\a u(x_0) \, x_\a$ and $\overline
b_0(x_\a):=\overline b(x_0)$. For every $\rho>0$, Theorem
\ref{bending} implies the existence of a sequence $\{v^\rho_n\}
\subset W^{1,p}(B'(x_0,\rho) \times I ;\Rb^3)$ such that $v^\rho_n
\to u_0$ in $L^p(B'(x_0,\rho) \times I;\Rb^3)$ (thus {\it a
fortiori} in $L^1(B'(x_0,\rho) \times I;\Rb^3)$),
$\frac{1}{\e_n}\int_I\nabla_3v^\rho_n(\cdot,x_3)\, dx_3 \wto
\overline b_0$ in $L^p(B'(x_0,\rho);\Rb^3)$ and
$$\lim_{n \to +\infty}\int_{B'(x_0,\rho) \times I} W_{\e_n}\left(x,\nabla_\a
v^\rho_n \Big|\frac{1}{\e_{n}} \nabla_3 v^\rho_n \right) dx =
\int_{B'(x_0,\rho)} W^*(x_\a,\nabla_\a u(x_0)|\overline b(x_0))\,
dx_\a.$$ Since $u_0 \in L^\infty(\o;\Rb^3)$, by Lemma \ref{bd} and
Remark \ref{sobo}, for any $\eta>0$ we can find a sequence
$\{w_n^\rho\} \subset W^{1,p}(B'(x_0,\rho) \times I ;\Rb^3)$ and
$C_\rho>0$ such that $\sup_n \|w_n^\rho\|_{L^\infty(B'(x_0,\rho)
\times I;\Rb^3)} \leq C_\rho$, $w^\rho_n \to u_0$ in
$L^p(B'(x_0,\rho) \times I;\Rb^3)$, $\frac{1}{\e_n}\int_I\nabla_3
w^\rho_n(\cdot,x_3)\, dx_3 \wto \overline b_0$ in
$L^p(B'(x_0,\rho);\Rb^3)$ and
\begin{eqnarray*}
&&\limsup_{n \to +\infty}\int_{B'(x_0,\rho) \times I}
W_{\e_n}\left(x, \nabla_\a w^\rho_n \Big|\frac{1}{\e_{n}} \nabla_3
w^\rho_n \right) dx\\
&&\hspace{4cm} \leq  \int_{B'(x_0,\rho)} W^*(x_\a,\nabla_\a
u(x_0)|\overline b(x_0))\, dx_\a + \LL^2(B'(x_0,\rho))
\eta.\end{eqnarray*} Thanks to (\ref{pgW}) and the separately convex
character of $W^*(x_0,\cdot|\cdot)$ (see the proof of Theorem
\ref{bending}), it follows that $W^*(x_0,\cdot|\cdot)$ is
$p$-Lipschitz. Thus our choice of $x_0$ implies that
\begin{equation}\label{ledret}\limsup_{\rho \to 0} \limsup_{n \to
+\infty}\med_{B'(x_0,\rho) \times I} W_{\e_n}\left(x, \nabla_\a
w^\rho_n \Big|\frac{1}{\e_{n}} \nabla_3 w^\rho_n \right) dx \leq
W^*(x_0,\nabla_\a u(x_0)|\overline b(x_0)) +\eta\end{equation} and
from the coercivity condition (\ref{pg1}), we get
\begin{equation}\label{coerc}
\sup_{\rho>0,\, n \in \Nb}\med_{B'(x_0,\rho) \times I} \left| \left(
\nabla_\a w^\rho_n \Big|\frac{1}{\e_{n}} \nabla_3 w^\rho_n
\right)\right|^p dx < +\infty.
\end{equation}
Let $\overline b_k \in \C^\infty_c(\o;\Rb^3)$ be such that
$\overline b_k \to \overline b$ in $L^p(\o;\Rb^3)$ and define
$$u_{n,k}^\rho(x):= u(x_\a) +\e_nx_3 (\overline b_k(x_\a) - \overline
b(x_0))+ w^\rho_n(x_\a,x_3) - \nabla_\a u(x_0) \, x_\a.$$ Then,
$u_{n,k}^\rho \to u$ in $L^1(B'(x_0,\rho) \times I;\Rb^3)$,
$\frac{1}{\e_n}\int_I \nabla_3 u^\rho_{n,k} (\cdot,x_3)\, dx_3 \wto
\overline b_k$ in $L^p(B'(x_0,\rho);\Rb^3)$ as $n \to +\infty$ and
$\sup_n \|u_{n,k}^\rho\|_{L^\infty(B'(x_0,\rho) \times I;\Rb^3)}
<+\infty$. Thus, since $S_{u_{n,k}^\rho} \cap (B'(x_0,\rho) \times
I) = (S_u \cap B'(x_0,\rho)) \times I$, we get that
\begin{eqnarray*}
\E_\infty(u,\overline b_k,B'(x_0,\rho)) & \leq & \liminf_{n \to
+\infty}\left\{ \int_{B'(x_0,\rho) \times I} W_{\e_n}\left(x,\nabla_\a u_{n,k}^\rho \Big|\frac{1}{\e_n} \nabla_3 u_{n,k}^\rho \right) dx \right.\\
&&\hspace{2.0cm} \left. + \int_{S_{u_{n,k}^\rho} \cap (B'(x_0,\rho)
\times I)} \left| \left( \left(\nu_{u_{n,k}^\rho}\right)_\a \Big|\frac{1}{\e_n} \left(\nu_{u_{n,k}^\rho}\right)_3 \right)\right| d\mathcal H^2 \right\}\\
& \leq & \liminf_{n \to +\infty} \int_{B'(x_0,\rho) \times I}
W_{\e_n}\Big(x, \nabla_\a u(x_\a) - \nabla_\a u(x_0) + \nabla_\a w^\rho_n (x)\\
&&\hspace{3cm} +\e_n x_3 \nabla_\a \overline b_k(x_\a)
\Big|\frac{1}{\e_{n}} \nabla_3 w^\rho_n(x) + \overline b_k(x_\a) -
\overline b(x_0)\Big)\, dx \\
&&\hspace{2.0cm} + \mathcal H^1(S_u \cap B'(x_0,\rho)).
\end{eqnarray*}
Thus from (\ref{sing}), we obtain
\begin{eqnarray*}
f(x_0) & \leq & \liminf_{\rho \to 0} \liminf_{k \to +\infty}
\liminf_{n \to +\infty}  \med_{B'(x_0,\rho) \times
I} W_{\e_n}\Big(x,\nabla_\a u(x_\a) - \nabla_\a u(x_0) + \nabla_\a w^\rho_n (x)\\
&&\hspace{3cm} +\e_n x_3 \nabla_\a \overline b_k(x_\a)
\Big|\frac{1}{\e_{n}} \nabla_3 w^\rho_n(x) + \overline b_k(x_\a) -
\overline b(x_0)\Big)\, dx. \end{eqnarray*} Relations
(\ref{pg1bis}), (\ref{ledret}), (\ref{coerc}) and H\"older's
inequality yield
\begin{eqnarray*}
f(x_0) & \leq & \liminf_{\rho \to 0} \liminf_{k \to +\infty}
\liminf_{n \to +\infty} \Bigg\{  \med_{B'(x_0,\rho) \times I}
W_{\e_n}\left(x,\nabla_\a w^\rho_n  \Big|\frac{1}{\e_n} \nabla_3 w^\rho_n \right) dx\\
&&\hspace{0cm} + c \med_{B'(x_0,\rho) \times I}
\bigg(1+|\nabla_\a u(x_\a) - \nabla_\a u(x_0)|^{p-1}+ |\overline b_k(x_\a) - \overline b(x_0)|^{p-1}\\
&&\hspace{0cm}   + \Big| \Big(  \nabla_\a w^\rho_n(x)
\Big|\frac{1}{\e_{n}} \nabla_3 w^\rho_n(x) \Big)\Big|^{p-1} +
\e_n^{p-1}|\nabla_\a \overline b_k(x_\a)|^{p-1}\bigg)\bigg(
|\nabla_\a u(x_\a) - \nabla_\a u(x_0)|\\
&&\hspace{1cm} +\e_n|\nabla_\a \overline b_k(x_\a)|+ |\overline b_k(x_\a) - \overline b(x_0)|\bigg)dx \Bigg\}\\
& \leq & W^*(x_0,\nabla_\a u(x_0)|\overline b(x_0)) +\eta\\
&&+ c\,  \limsup_{\rho \to 0}  \left\{ \med_{B'(x_0,\rho)}
\big(1+ |\nabla_\a u(x_\a) - \nabla_\a u(x_0)|^p + |\overline b(x_\a)-\overline b(x_0)|^p \big)\, dx_\a \right\}^{(p-1)/p}\\
&&\hspace{2cm} \times \left\{\med_{B'(x_0,\rho)} \big( |\nabla_\a
u(x_\a) - \nabla_\a u(x_0)|^p + |\overline b(x_\a)-\overline
b(x_0)|^p \big)\, dx_\a \right\}^{1/p}.
\end{eqnarray*}
Thanks to our choice of $x_0$ and letting $\eta \to 0$, we conclude
that $f(x_0) \leq W^*(x_0,\nabla_\a u(x_0)|\overline b(x_0))$ for
$\mathcal L^2$-a.e. $x_0 \in \o$ which completes the proof in the
case where $u \in L^\infty(\o;\Rb^3) \cap SBV^p(\o;\Rb^3)$. The
general case can in turn be treated by approximation exactly as in
the proof of \cite[Lemma~3.8]{babadjian}.
\end{proof}

\subsection{The lower bound}\label{lower}

Let us now prove the lower bound. The proof is essentially based on
Theorem \ref{bab} and a blow up argument.

\begin{lemma}\label{gi}
For all $u \in BV(\O;\Rb^3)$ and all $\overline b \in
L^p(\o;\Rb^3)$, $\E(u,\overline b,\o) \geq \mathcal G(u,\overline
b)$.
\end{lemma}

\begin{proof}It is not restrictive to assume that $\E(u,\overline
b,\o)<+\infty$. By $\G$-convergence, there exists a sequence
$\{u_n\} \subset SBV^p(\O;\Rb^3)$ such that $u_n \to u$ in
$L^1(\O;\Rb^3)$, $\frac{1}{\e_n}\int_I \nabla_3u_n(\cdot,x_3)\, dx_3
\wto \overline b$ in $L^p(\o;\Rb^3)$ and
\begin{equation}\label{1805}
\lim_{n \to +\infty}\left[ \int_{\O} W_{\e_n}\left(x, \nabla_\a u_n
\Big|\frac{1}{\e_n}\nabla_3 u_n \right)dx + \int_{S_{u_n}}\left|
\left( \big( \nu_{u_n} \big)_\a \Big|\frac{1}{\e_n} \big( \nu_{u_n}
\big)_3 \right) \right|d\mathcal H^2 \right] = \E(u,\overline b,\o).
\end{equation}
Arguing exactly as in the proof of \cite[Lemma~3.9]{babadjian}, we
can actually show that $u \in SBV^p(\o;\Rb^3)$ and that $u_n \wto u$
in $SBV^p(\O;\Rb^3)$. Now for every Borel set $E \subset \o$, define
the following sequences of Radon measures:
$$\lambda_n(E):=W_{\e_n}\left(\cdot, \nabla_\a u_n
\Big|\frac{1}{\e_n}\nabla_3 u_n \right)\LL^3\res \,  (E \times I) +
\left| \left( \big( \nu_{u_n} \big)_\a \Big|\frac{1}{\e_n} \big(
\nu_{u_n} \big)_3 \right) \right|\mathcal H^2\res \,  (S_{u_n}
\cap(E \times I))$$ and
$$\mu_n(E):=\left| \left( \big( \nu_{u_n} \big)_\a \Big|\frac{1}{\e_n} \big(
\nu_{u_n} \big)_3 \right) \right|\mathcal H^2\res \,  (S_{u_n}
\cap(E \times I)).$$ Then for a subsequence (not relabeled), there
exist nonnegative and finite Radon measures $\lambda$ and $\mu \in
\M(\o)$ such that $\lambda_n \xrightharpoonup[]{*} \lambda$ and
$\mu_n \xrightharpoonup[]{*} \mu$ in $\M(\o)$. By the Besicovitch
Differentiation Theorem (\cite[Theorem~2.22]{AFP}), one can find
three mutually disjoint nonnegative Radon measures $\lambda^a$,
$\lambda^j$ and $\lambda^c$ such that
$\lambda=\lambda^a+\lambda^j+\lambda^c$ where $\lambda^a \ll \LL^2$
and $\lambda^j \ll \HH^1\res \,  S_u$. It is enough to check that
\begin{equation}\label{surfinf}
\frac{d\lambda^j}{d\HH^1\res \,  S_u}(x_0) \geq 1 ,\quad \text{ for
}\HH^1\text{-a.e. }x_0 \in S_u\end{equation} and
\begin{equation}\label{bulkinf}
\frac{d\lambda^a}{d\LL^2}(x_0) \geq  W^*(x_0,\nabla_\a
u(x_0)|\overline b(x_0)),\quad \text{ for }\LL^2\text{-a.e. }x_0 \in
\o.\end{equation} Indeed, if (\ref{surfinf}) and (\ref{bulkinf})
hold, we obtain from (\ref{1805}) that
\begin{eqnarray*}
\E(u,\overline b,\o) & \geq & \lambda(\o) =
\lambda^a(\o)+\lambda^j(\o) + \lambda^c(\o)\\
& \geq &  \int_\o W^*(x_\a,\nabla_\a u|\overline b)\, dx_\a +
\HH^1(S_u) = \mathcal G(u,\overline b).
\end{eqnarray*}

We first prove (\ref{surfinf}). Fix a point $x_0 \in S_u$ such that
$$\frac{d\lambda^j}{d\HH^1\res \,  S_u}(x_0)=\frac{d\lambda}{d\HH^1\res \,  S_u}(x_0)$$
exists and is finite and remark that $\HH^1$-a.e. points in $S_u$
satisfy this property. Let $\{\rho_k\} \searrow 0^+$ be such that
$\lambda(\partial B'(x_0,\rho_k))=0$ for each $k\in \Nb$. Then,
\begin{eqnarray*}
\frac{d\lambda}{d\HH^1\res \,  S_u}(x_0) & = & \lim_{k \to
+\infty}\frac{\lambda(B'(x_0,\rho_k))}{\HH^1(S_u
\cap B'(x_0,\rho_k))}\\
& = & \lim_{k \to +\infty}\lim_{n \to
+\infty}\frac{\lambda_n(B'(x_0,\rho_k))}{\HH^1(S_u \cap
B'(x_0,\rho_k))}\\
& \geq & \liminf_{k \to +\infty}\liminf_{n \to
+\infty}\frac{\HH^2(S_{u_n} \cap (B'(x_0,\rho_k) \times
I))}{\HH^1(S_u \cap B'(x_0,\rho_k))}.
\end{eqnarray*}
By \cite[Theorem~4.36]{AFP}, we have that $$\liminf_{n \to
+\infty}\HH^2(S_{u_n} \cap (B'(x_0,\rho_k) \times I)) \geq \HH^1(S_u
\cap
B'(x_0,\rho_k))$$ hence we obtain (\ref{surfinf}).\\

Let us prove that (\ref{bulkinf}) holds at every point $x_0 \in \o
\setminus N$ (where $N \subset \o$ is the exceptional set introduced
in Proposition \ref{rho}) which is a Lebesgue point of both
$\nabla_\a u$ and $\overline b$, a point of approximate
differentiability of $u$ such that
$$\frac{d\lambda^a}{d\LL^2}(x_0)=\frac{d\lambda}{d\LL^2}(x_0)$$
exists and is finite and satisfying
\begin{equation}\label{choice}\lim_{\rho \to
0}\frac{\mu(B'(x_0,\rho))}{2\rho}=0.\end{equation} It turns out that
$\LL^2$-a.e. points $x_0$ in $\o$ satisfy these property. Indeed,
the verification of (\ref{choice}) is similar to the one of
(\ref{itemc}) used in the proof of Theorem \ref{Ambrosio}. As
before, let $\{\rho_k\} \searrow 0^+$ be such that $\lambda(\partial
B'(x_0,\rho_k))=0$ for every $k \in \Nb$, then
\begin{eqnarray}\label{abscontpart}
\frac{d\lambda}{d\LL^2}(x_0)  & = &
\lim_{k \to +\infty}\frac{\lambda(B'(x_0,\rho_k))}{\LL^2(B'(x_0,\rho_k))}\nonumber\\
& = & \lim_{k \to +\infty}\lim_{n \to
+\infty}\frac{\lambda_n(B'(x_0,\rho_k))}{\LL^2(B'(x_0,\rho_k))}\nonumber\\
& \geq & \limsup_{k \to +\infty}\limsup_{n \to
+\infty}\frac{1}{\LL^2(B'(x_0,\rho_k))}\int_{B'(x_0,\rho_k) \times
I}W_{\e_n}\left(x,\nabla_\a u_n \Big|\frac{1}{\e_n}\nabla_3 u_n \right)dx\nonumber\\
& = & \limsup_{k \to +\infty}\limsup_{n \to
+\infty}\frac{1}{\LL^2(B')}\int_{B' \times
I}W_{\e_n}\left(x_0+\rho_kx_\a,x_3,\nabla_\a u_{n,k}
\Big|\frac{\rho_k}{\e_n}\nabla_3 u_{n,k} \right)dx,
\end{eqnarray}
where $u_{n,k}(x_\a,x_3)=[u_n(x_0+\rho_kx_\a,x_3) - u(x_0)]/\rho_k$.
Since $x_0$ is a point of approximate differentiability of $u$, we
have that
\begin{equation}\label{approxdiff}
\lim_{k \to +\infty}\lim_{n \to +\infty} \int_{B' \times
I}|u_{n,k}(x) - \nabla_\a u(x_0)x_\a|\, dx =0
\end{equation}
and using the fact that $x_0$ is a Lebesgue point of $\overline b$,
for every $v \in L^{p'}(B';\Rb^3)$ we get that
\begin{equation}\label{weakconvb}
\lim_{k \to +\infty}\lim_{n \to +\infty} \int_{B'}\left( \frac{
\rho_k}{\e_n}\int_I \nabla_3 u_{n,k}(x_\a,x_3)\, dx_3\right)\cdot
v(x_\a)\, dx_\a= \int_{B'}\overline b(x_0)\cdot v\, dx_\a.
\end{equation}
Changing variables in the surface term and thanks to (\ref{choice}),
it yields to
\begin{eqnarray}\label{lastsurf}
&&\limsup_{k \to +\infty}\limsup_{n \to +\infty} \int_{S_{u_{n,k}}
\cap (B' \times I)}\left| \left( \big( \nu_{u_{n,k}} \big)_\a
\Big|\frac{\rho_k}{\e_n} \big( \nu_{u_{n,k}} \big)_3 \right)
\right|d\mathcal H^2\nonumber\\
&&\hspace{2cm} = \limsup_{k \to +\infty}\limsup_{n \to +\infty}
\frac{1}{\rho_k} \int_{S_{u_n} \cap (B'(x_0,\rho_k) \times I)}\left|
\left( \big( \nu_{u_n} \big)_\a \Big|\frac{1}{\e_n} \big( \nu_{u_n}
\big)_3 \right) \right|d\mathcal H^2\nonumber\\
&&\hspace{2cm} \leq \limsup_{k \to +\infty}\limsup_{n \to
+\infty}\frac{\mu_n(\overline{B'(x_0,\rho_k)})}{\rho_k}\nonumber\\
&&\hspace{2cm} \leq \limsup_{k \to +\infty}
\frac{\mu(\overline{B'(x_0,\rho_k)})}{\rho_k}=0
\end{eqnarray}
because $\mu(\partial B'(x_0,\rho_k)) \leq \lambda(\partial
B'(x_0,\rho_k))=0$. Set
\begin{equation}\label{borneM}
M:=\max\left\{
\left(\frac{\LL^2(B')}{\b'}\left(\left|\frac{d\lambda}{d\LL^2}(x_0)\right|+1\right)\right)^{1/p},
|\overline b(x_0)|\LL^2(B')^{1/p}\right\} <+\infty.
\end{equation}
From (\ref{abscontpart})-(\ref{lastsurf}), using a diagonalization
argument, the fact that $L^{p'}(B';\Rb^3)$ is separable and Remark
\ref{diag}, we can find a sequence $n(k) \nearrow +\infty$ such
that, setting $\d_k:=\e_{n(k)}/\rho_k$, $v_k:=u_{n(k),k}$,
$u_0(x_\a):=\nabla u(x_0)\, x_\a$ and $\overline
b_0(x_\a):=\overline b(x_0)$, then $\d_k\to 0$, $v_k \to u_0$ in
$L^1(B'\times I;\Rb^3)$, $\frac{1}{\d_k}\int_I \nabla_3
v_k(\cdot,x_3)\, dx_3 \wto \overline b_0$ in $L^p(B';\Rb^3)$,
\begin{equation}\label{1557}
\lim_{k \to +\infty} \int_{S_{v_k}}\left| \left( \big( \nu_{v_k}
\big)_\a \Big|\frac{1}{\d_k} \big( \nu_{v_k} \big)_3 \right)
\right|d\mathcal H^2=0,\end{equation}
\begin{equation}\label{1558}\frac{d\lambda}{d\LL^2}(x_0) \geq \limsup_{k \to +\infty}\frac{1}{\LL^2(B')}
\int_{B' \times I}W_{\e_{n(k)}}\left(x_0+\rho_k x_\a,x_3,\nabla_\a
v_k \Big|\frac{1}{\d_k}\nabla_3 v_k \right)dx\end{equation} and for
every $(u,\overline b) \in L^p(B' \times I;\Rb^3)\times
L^p(B';\Rb^3)$ with $\|\overline b\|_{L^p(B';\Rb^3)} \leq M$, the
$\G$-limit in $L^p_s(B' \times I;\Rb^3) \times L^p_w(B';\Rb^3)$ of
$$
\left\{\begin{array}{ll} \ds \int_{B'\times I}
W_{\e_{n(k)}}\left(x_0+\rho_k x_\a,x_3,\nabla_\a
u\Big|\frac{1}{\d_k}\nabla_3 u \right)\, dx & \text{ if } \left\{
\begin{array}{l}
u \in W^{1,p}(B'\times I;\Rb^3),\\
\overline b=\frac{1}{\d_k}\int_I \nabla_3 u(\cdot,x_3)\, dx_3,
\end{array}
\right.\\[0.3cm]
+\infty & \text{ otherwise},
\end{array}\right.$$
coincides with
$$ \left\{\begin{array}{ll} \ds \int_{B'} W^*(x_0,\nabla_\a u|\overline
b)\, dx_\a & \text{ if } u \in W^{1,p}(B';\Rb^3),\\[0.3cm]
+\infty & \text{ otherwise}.
\end{array}\right.$$
From (\ref{1557}), (\ref{1558}) and (a slight variant of) Lemma
\ref{bd}, for any $0<\eta<1$, there exist a constant $C>0$ and
$\{w_k\} \subset SBV^p(B' \times I;\Rb^3)$ such that $w_k \to u_0$
in $L^1(B'\times I;\Rb^3)$, $\frac{1}{\d_k}\int_I \nabla_3
w_k(\cdot,x_3)\, dx_3 \wto \overline b_0$ in $L^p(B';\Rb^3)$,
$\sup_k \|w_k\|_{L^\infty(B'\times I;\Rb^3)} \leq C$,
$$\lim_{k \to +\infty}
\int_{S_{w_k}}\left| \left( \big( \nu_{w_k} \big)_\a
\Big|\frac{1}{\d_k} \big( \nu_{w_k} \big)_3 \right) \right|d\mathcal
H^2=0$$ and $$\frac{d\lambda}{d\LL^2}(x_0)  \geq \limsup_{k \to
+\infty}\frac{1}{\LL^2(B')}\int_{B' \times
I}W_{\e_{n(k)}}\left(x_0+\rho_k x_\a, x_3, \nabla_\a w_k
\Big|\frac{1}{\d_k}\nabla_3 w_k \right)dx-\eta.$$ From the
$p$-coercivity condition (\ref{pg1}) and \cite[Theorem 4.36]{AFP},
the sequence $\{w_k\}$ converges weakly to $u$ in $SBV^p(\O;\Rb^3)$
and it fulfills the assumptions of Theorem \ref{bab}. Thus, for a
not relabeled subsequence, one can find another sequence $\{z_k\}
\subset W^{1,\infty}(B' \times I;\Rb^3)$ such that $z_k \wto u_0$ in
$W^{1,p}(B'\times I;\Rb^3)$, $\frac{1}{\d_k}\int_I \nabla_3
z_k(\cdot,x_3)\, dx_3 \wto \overline b_0$ in $L^p(B';\Rb^3)$,
$\big\{\big|\big( \nabla_\a z_k|\frac{1}{\d_k}\nabla_3 z_k\big)
\big|^p \big\}$ is equi-integrable and $\LL^3(\{z_k \neq w_k\} \cup
\{\nabla z_k \neq \nabla w_k\}) \to 0$. Hence
$$\frac{d\lambda}{d\LL^2}(x_0)  \geq \limsup_{k \to
+\infty}\frac{1}{\LL^2(B')}\int_{\{w_k=z_k\}}W_{\e_{n(k)}}\left(x_0+\rho_k
x_\a, x_3,\nabla_\a z_k \Big|\frac{1}{\d_k}\nabla_3 z_k
\right)dx-\eta$$ and using the $p$-growth condition (\ref{pg1}), the
fact that $\big\{\big|\big( \nabla_\a z_k|\frac{1}{\d_k}\nabla_3
z_k\big) \big|^p \big\}$ is equi-integrable and that $\LL^3(\{z_k
\neq w_k\} ) \to 0$ we get,
$$\limsup_{k \to +\infty} \int_{\{w_k\neq z_k\}}W_{\e_{n(k)}}\left(x_0+\rho_k x_\a, x_3,\nabla_\a
z_k \Big|\frac{1}{\d_k}\nabla_3 z_k \right)dx=0.$$ As a consequence
$$\frac{d\lambda}{d\LL^2}(x_0) \geq \limsup_{k \to
+\infty}\frac{1}{\LL^2(B')}\int_{B' \times
I}W_{\e_{n(k)}}\left(x_0+\rho_k x_\a, x_3,\nabla_\a z_k
\Big|\frac{1}{\d_k}\nabla_3 z_k \right)dx-\eta$$ and by the
$p$-coercivity condition (\ref{pg1}) and (\ref{borneM}),
$$\left\| \frac{1}{\d_k} \int_I \nabla_3 z_k(\cdot,x_3)\, dx_3 \right\|_{L^p(B';\Rb^3)} \leq M,
\quad \|\overline b_0\|_{L^p(B';\Rb^3)} \leq M.$$ Thus by our choice
of the subsequence $n(k)$ and Remark \ref{diag}, we get that
$$\frac{d\lambda}{d\LL^2}(x_0) \geq W^*(x_0,\nabla_\a u(x_0)|\overline b(x_0)) - \eta.$$
Letting $\eta$ tend to zero completes the proof of (\ref{bulkinf}).
\end{proof}

\begin{rmk}{\rm
Note that it seems difficult to think of applying the decoupling
variable method introduced in \cite{babF} and further developed in
\cite{BB1,BB2}. Indeed, this generalized framework has the drawback
that we have no information on the way that $W_\e$ depends on $\e$,
and it requires application of such abstract results as
metrizability of $\G$-convergence. Remark also that the same kind of
blow-up argument considered here could have been used in
\cite{BB1,BB2,babF} in place of the decoupling variable method, in
order to treat the presence of the spatial variable. }\end{rmk}

\section{case without bending moment}\label{7}

\noindent In this last section, we deduce from Theorem
\ref{gammaconvbend} a similar result without the presence of the
bending moment. Define $\I_\e : L^p(\O;\Rb^3) \to [0,+\infty]$ by
$$\I_\e(u):=\left\{\begin{array}{ll}
\ds \int_{\O} W_\e\left(x,\nabla_\a u \Big|\frac{1}{\e}\nabla_3
u\right)dx & \text{ if }u \in W^{1,p}(\O;\Rb^3),\\[0.3cm]
 +\infty & \text{ otherwise.} \end{array}\right.$$
In \cite[theorem~2.5]{BFF}, it has been proved the following
integral representation result:
\begin{thm}\label{nobending}
For every sequence $\{\e_n\} \searrow 0^+$, there exist a
subsequence (not relabeled) and a Carath\'eodory function $\widehat
W : \o \times \Rb^{3 \times 2} \to [0,+\infty)$ (depending on the
subsequence) such that the sequence $\I_{\e_n}$ $\G$-converges in
$L_s^p(\O;\Rb^3)$ to $\I$ where
$$\I(u)=\left\{
\begin{array}{ll}
\ds \int_\o \widehat W(x_\a,\nabla_\a u)\, dx_\a & \text{ if
}u \in W^{1,p}(\o;\Rb^3),\\[0.3cm]
+\infty & \text{ otherwise}.
\end{array}
\right.$$
\end{thm}

We refer to \cite{LDR,BFF,babF,BB1,BB2} for more explicit formulas
in particular cases.

\begin{rmk}\label{min}{\rm
As it has been pointed out in \cite{babF} in the case where $W_\e$
was independent of $\e$ (see also \cite{BFMbis}), it can still be
seen here that
$$\widehat W(x_0,\overline \xi)=\min_{z \in \Rb^3}W^*(x_0,\overline
\xi|z)$$ for all $\overline \xi \in \Rb^{3 \times 2}$ and a.e. $x_0
\in \o$.}
\end{rmk}

Define now $\F_\e : BV(\O;\Rb^3) \to [0,+\infty]$ by
$$\F_\e(u):=\left\{
\begin{array}{ll}
\ds \int_\O W_\e\left(x,\nabla_\a u \Big|\frac{1}{\e}\nabla_3
u\right)dx + \int_{S_u} \left| \left( \left(\nu_u\right)_\a
\Big|\frac{1}{\e} \left( \nu_u \right)_3 \right)\right| d\HH^2 & \text{if }u \in SBV^p(\O;\Rb^3),\\[0.3cm]
\ds +\infty & \text{otherwise}.
\end{array}
\right.$$

As a consequence of Theorem \ref{gammaconvbend}, Theorem
\ref{nobending}, Remark \ref{min} and a standard measurability
selection criterion (see {\it e.g.}
\cite[Theorem~1.2,~Chapter~VIII]{ET}) we get the following integral
representation result for dimension reduction problems in $SBV$
without bending moment:

\begin{thm}\label{gammaconv}
For every sequence $\{\e_n\} \searrow 0^+$, there exists a
subsequence, still labeled $\{\e_n\}$ such that $\F_{\e_n}$
$\G$-converges in $L_s^1(\O;\Rb^3)$ to $\F : BV(\O;\Rb^3)\to
[0,+\infty]$ defined by
$$\F(u):=\left\{
\begin{array}{ll}
\ds \int_\o \widehat W(x_\a,\nabla_\a u)\, dx_\a +\HH^1(S_u) & \text{if }  u \in SBV^p(\o;\Rb^3),\\[0.3cm]
\ds  +\infty & \text{otherwise},
\end{array}
\right.$$ where $\widehat W$ is given by Theorem \ref{nobending}.
\end{thm}

\vspace{0.5cm}

\noindent {\it Acknowledgments.} The author wishes to thank Irene
Fonseca and Gilles Francfort for having drawn this problem to his
attention. The research of J.-F. Babadjian has been supported by the
MULTIMAT Marie Curie Research Training Network MRTN-CT-2004-505226
``Multi-scale modelling and characterisation for phase
transformations in advanced materials''.

\end{document}